\documentclass[11pt,a4paper,twoside]{amsart}

\usepackage{amsmath,amssymb,amsfonts,amsthm,mathrsfs}
\usepackage{times,hyperref,color}
\usepackage{algorithm,algorithmic}

\newcommand{\K}{\mathbb{K}}

\let\mathcal\mathscr

\newtheorem{The}{Theorem}[section]

\newtheorem{Theorem}{Theorem}[section]
\newtheorem{Proposition}[The]{Proposition}

\theoremstyle{definition}
\newtheorem{Definition}[The]{Definition}

\newtheorem{Example}[The]{Example}


\begin{document}

\title{
A Gr\"obner-bases algorithm for the computation
\\
of the cohomology of Lie (super) algebras}

\author{Mansour Aghasi}
\address{Department of Mathematical Sciences,
Isfahan University of Technology, Isfahan, IRAN}
\email{m.aghasi@cc.iut.ac.ir}

\author{Benyamin M.-Alizadeh}
\address{Department of Mathematical Sciences,
Isfahan University of Technology, Isfahan, IRAN} 
\email{b.alizadeh@math.iut.ac.ir}

\author{Jo\"{e}l Merker}
\address{D\'epartment de Math\'ematiques d'Orsay,
B\^atiment 425, Facult\'e des Sciences, Universit\'e Paris XI - Orsay, 
F-91405 Orsay Cedex, FRANCE}
\email{merker@dma.ens.fr}

\author{Masoud Sabzevari}
\address{Department of Mathematical Sciences,
Isfahan University of Technology, Isfahan, IRAN}
\email{sabzevari@math.iut.ac.ir}

\date{\number\year-\number\month-\number\day}

\maketitle

\begin{abstract}
We present an effective algorithm for computing the standard
cohomology spaces of finitely generated Lie (super) algebras over
a commutative field $\K$ of characteristic zero. 
In order
to reach explicit representatives of some generators of the quotient
space $\mathcal{ Z}^k \big / \mathcal{ B}^k$ of cocycles $\mathcal{
Z}^k$ modulo coboundaries $\mathcal{ B}^k$, we apply
Gr\"obner bases techniques (in the appropriate linear
setting) and take advantage of 
their strength. 
Moreover, when the considered Lie (super) algebras enjoy a
grading\,\,---\,\,a case which often happens both in representation
theory and in differential geometry\,\,---, all cohomology spaces
$\mathcal{ Z}^k \big / \mathcal{ B}^k$ naturally split up as direct sums
of smaller subspaces, and this enables us, for higher dimensional Lie
(super) algebras, to improve the computer speed of
calculations. Lastly, we implement our algorithm in the {\sc Maple}
software and evaluate its performances via some examples, most of
which have several applications in the theory of Cartan-Tanaka
connections.
\end{abstract}

\pagestyle{headings} 
\markleft{\sf Mansour Aghasi, Benyamin M.-Alizadeh, Jo\"el Merker
and Masoud Sabzevari}
\markright{\sf A Gr\"obner-bases algorithm for the computation
of the cohomology of Lie (super) algebras}

\section{Introduction}
\label{introduction}

The concept of cohomology group\,\,---\,\,one of the
central concepts in contemporary science\,\,---\,\,possesses
established applications in several areas of pure mathematics, for
instance: deformation of Lie algebras (\cite{ Goze}); analytic
partial differential equations; global foliation theory; combinatorics
(Mcdonald identities); invariant differential operators; cobordism
theory; infinite-dimensional Lie algebras (\cite{Fuks}); exterior
differential systems; Cartan-Tanaka theory of connections (\cite{BES,
AMS, MS}); {\em etc.} Moreover, cohomology groups also have
applications in
quantum physics; for quasi-invariancy of certain Lagrangians; in the
Wess-Zumino-Novikov-Witten
model ({\em cf.}~\cite{Azcarraga}); when one reinterprets
general relativity by means of $\mathfrak{ so}(3,1)$-valued
connections; {\em etc.} It therefore turns out to be worthwhile to set
up appropriate efficient algorithms for the computation of Lie (super)
algebra cohomologies, granted that calculations quickly become hard by
hand.

Recently, a few articles have been published in this
direction. Kornyak \cite{Kornyak2002,Kornyak2004} devised an algorithm
and implemented it in the ${\sf C}$ 
program. Moreover, Grozman, Leites, Post
and Von Hijligenberg (\cite{Grozman, Leites, Post}) prepared some
packages for computing Lie (super) algebra cohomologies in {\sc
Reduce} and in {\sc Mathematica}. In the present article, motivated by
the specific objective of developing the construction of {\em
effective} Cartan-Tanaka connections that are valued in Lie algebras
which are {\em not} semi-simple ({\em see} \cite{BES, AMS, MS} for some
instances of that research program and also~\cite{ CapSlovak} in
the parabolic$/$simple case), our main aim is to set up an alternative
algorithm and to implement it in the {\sc Maple} software. We would
like to employ the method of {\em Gr\"obner bases}, a modern,
effective and widespread tool in computational mathematics. Of course,
the continued regular progresses in Gr\"obner bases
algorithms enrich {\em
de facto} any algorithm that is built on them. For convenience and
self-contentness, a short reminder of Gr\"obner bases concepts will be
given in Section~2. But before that, let us present a brief
description of the definitions, notations and formulas in Lie super
algebras, and let us introduce their cohomology groups, precisely.

A {\sl Lie super algebra}
over a commutative field $\mathbb K$ of characteristic zero
is a $(\mathbb{ Z} / 2\mathbb{Z} )$-graded algebra which is
a direct sum (as a vector space):
\[
\frak g
=
\frak g_{\overline{0}}
\oplus
\frak g_{\overline{1}}
\]
of two subspaces $\mathfrak{ g}_{ 
\overline{ 0}}$ and $\mathfrak{ g}_{ 
\overline{ 1}}$, together with a {\sl degree-zero}
graded Lie bracket: 
\[
[\cdot,\cdot]\,
\colon\ \
\frak g\times\frak g\longrightarrow\frak g, 
\] 
that is to say: $[ \cdot , \cdot ]$ is a bilinear map satisfying:
\[
\big[
\frak g_{\overline{i}},\frak g_{\overline{j}}
\big]
\subseteq
\frak g_{\overline{i+j}},
\] 
for any $i,j=0,1$ where $\overline{i+j}=i+j \,\, {\rm mod} \ 2$, and
satisfying also, for arbitrary elements ${\sf x},{\sf y},{\sf
z}\in\frak g$, the two standard conditions:
\[
[{\sf x},{\sf y}]=-(-1)^{|{\sf x}||{\sf y}|}[{\sf y},{\sf x}] 
\ \ \ \ \ \ \ \ \ \ \ \ \ \ \ \ \ \ \ \ \ \ \ \ \ 
\text{\footnotesize\sf (skew-symmetry)},
\]
\[
\big[{\sf x},[{\sf y},{\sf z}]\big]
= 
\big[
[{\sf x},{\sf y}],{\sf z}
\big]
+
(-1)^{|{\sf x}||{\sf y}|}
\big[{\sf y},[{\sf x},{\sf z}]\big] \ \
\ \ \ \ \ \ \ 
\text{\footnotesize\sf (Jacobi identity)},
\]
where the weight $|{\sf x}|$ is defined to be $0$ when ${\sf
x}\in\frak g_{\overline{0}}$ and to be $1$ when ${\sf x}\in\frak
g_{\overline{1}}$. The elements of $\frak g_{\overline{0}}$ and of
$\frak g_{\overline{1}}$ are called {\sl even} and {\sl odd},
respectively. In differentialo-geometric applications (\cite{ BES,
CapSlovak, AMS, MS}), the commutative
field $\mathbb{ K}$ of characteristic zero
is usually assumed to be
either just $\mathbb{ Q}$, or $\mathbb R$, or $\mathbb C$, plainly.

A {\sl $\mathfrak{ g}$-module} $V$ is a vector space over the same
field $\mathbb K$ together with a bilinear map (denoted shortly with a
dot) $\cdot \colon \frak g\times V\rightarrow V$ having the property:
\[
[{\sf x},{\sf y}]\cdot v
=
{\sf x}\cdot 
({\sf y}\cdot v)-(-1)^{|{\sf x}||{\sf y}|}{\sf y}\cdot 
({\sf x}\cdot v),
\]
for any two ${\sf x},{\sf y}\in\frak g$ and any $v\in V$. One of the
most important instances of such $\mathfrak{ g}$-modules occurs when
$\frak g$ happens to be a Lie (super) subalgebra of a certain larger
Lie (super) algebra $\frak h =: V$, with the bilinear map $\cdot
\colon \frak
g\times\frak h\rightarrow\frak h$ being just precisely the Lie bracket
of $\frak h$, of course.

Thus, let $\frak g$ be an $m$-dimensional Lie super algebra and let
$V$ be a $\frak g$-module. For any integer $k \geqslant 0$, the space
$\mathcal{ C}^k ( \mathfrak{ g},V)$ of {\sl $k$-cochains} consists of
the space of $k$-linear {\sl super skew-symmetric} maps:
\[
\Phi\colon\ \
\mathfrak g^k\longrightarrow V, 
\]
where $\mathfrak{ g}^k = \mathfrak{ g} \times \cdots \times \mathfrak{
g}$ ($k$ times, with $\mathfrak{ g}^0 = \{ 0\}$ naturally), and where
super skew-symmetry means symmetry with respect to the transposition
of odd elements and usual skew-symmetry with respect to all other
transpositions, that is to say generally:
\[
\Phi\big({\sf z}_1,\ldots,
{\sf z}_i,{\sf z}_{i+1},\ldots,
{\sf z}_k\big)
=
-(-1)^{|{\sf z}_i||{\sf z}_{i+1}|}\,
\Phi\big({\sf z}_1,\ldots,{\sf z}_{i+1},
{\sf z}_{i},\ldots,{\sf z}_k\big).
\]
Then for any integer $k \geqslant 0$, there
is a fundamental linear {\sl differential operator:}
\[
\partial^k\colon\ \ \
\mathcal{C}^k\big(\mathfrak{g},\,V\big) \longrightarrow
\mathcal{C}^{k+1}\big(\mathfrak{g},\,V\big),
\]
mapping a $k$-cochain $\Phi$ uniquely 
to a $(k+1)$-cochain $\partial^k \Phi$
that acts as follows ({\it see} \cite{Fuks, Khakim}) on any collection
of $k+1$ elements ${\sf e}_0, \dots, {\sf e}_p\in\frak
g_{\overline{0}}$, and ${\sf o}_{p+1}, \dots, {\sf o}_k\in\frak
g_{\overline{1}}$:
\begin{equation}
\label{cohomology} 
\small
\aligned
(\partial^k\Phi)
& 
\big({\sf e}_0,\dots,{\sf e}_p,\,
{\sf o}_{p+1}, \dots, {\sf o}_k\big) 
:= 
\\
&
:=
\sum_{i=0}^p\,(-1)^{i+1}\,
{\sf e}_i\cdot\Phi\big({\sf
e}_0,\dots,\widehat{{\sf e}}_i,\ldots, {\sf e}_p,\,{\sf o}_{p+1},\dots, 
{\sf o}_k\big) 
+
\\
&
+
\sum_{0\leqslant i<j\leqslant k}\, 
(-1)^{i+j+1}\, 
\Phi\big([{\sf e}_i,{\sf e}_j],{\sf e}_0, \dots,
\widehat{{\sf e}}_i,\ldots,
\widehat{{\sf e}}_j,\ldots,
{\sf e}_p,\,{\sf o}_{p+1}, \dots,{\sf o}_k\big)
+
\\
&
+
\sum_{i=0}^p\,\sum_{j=p+1}^k(-1)^i\,
\Phi\big({\sf e}_0,\dots,\widehat{{\sf e}}_i,\dots, 
{\sf e}_p,[{\sf e}_i,
{\sf o}_j],\,{\sf o}_{p+1},\dots,\widehat{{\sf o}}_j,\dots,
{\sf o}_k\big)
+
\\
&
+
\sum_{{p+1\leqslant i<j\leqslant k}}\,
\Phi\big([{\sf o}_i,{\sf o}_j],{\sf e}_0, \dots, {\sf
e}_p,\,{\sf o}_{p+1},\dots,
\widehat{{\sf o}}_i, \dots,\widehat{{\sf o}}_j,\ldots, {\sf o}_k\big)
+
\\
&+(-1)^p
\sum_{i=p+1}^{k}\,
{\sf o}_i\cdot
\Phi\big({\sf e}_0,\dots,\ldots,
{\sf e}_p,\,{\sf o}_{p+1},\dots,
\widehat{{\sf o}}_i,\ldots,{\sf o}_k\big),
\endaligned
\end{equation}
where as usual, $\widehat{\sf{z}}_l$ means removal of the term
${\sf{z}}_l$ (in the case of Lie algebras, comparing with some
references such as \cite{AMS,Azcarraga,Goze,MS}, there is an overall
minus sign in the right-hand side). One checks (\cite{Fuks}) that in
the case of Lie algebras $\mathfrak{ g} \subset \mathfrak{ h} = V$,
only the first two lines of the above definition are non-zero, and in
fact, for any $k+1$ vectors ${\sf z}_0, {\sf z}_1, \dots, {\sf z}_k
\in \mathfrak{ g}$, one has:
\begin{equation}
\label{standard-cochain}
\small
\aligned
(\partial^k\Phi)
\big({\sf z}_0,{\sf z}_1,\dots,{\sf z}_k\big)
&
:=
\sum_{i=0}^k\,(-1)^i
\big[{\sf z}_i,\,\Phi({\sf z}_0,
\dots,
\widehat{\sf z}_i,\dots,{\sf z}_k)\big]
+
\\
&
\ \ \ \ \
+
\sum_{0\leqslant i<j\leqslant k}\,
(-1)^{i+j}\,
\Phi\big([{\sf z}_i,{\sf z}_j],
{\sf z}_0,\dots,\widehat{\sf z}_i,
\dots,
\widehat{\sf z}_j,\dots,{\sf z}_k\big).
\endaligned
\end{equation}
In both cases, 
this $(k+1)$-cochain $\partial^k \Phi$ is
clearly linear with respect to each argument, and furthermore, it is
(super) skew-symmetric (\cite{ Fuks}). Furthermore, one can verify that
the compositions $\partial^{ k+1} \circ \partial^k$ vanish for any
$k\in\mathbb N$, hence we have the following {\sl cochain complex}:
\begin{equation}
\label{complex-partial}
0 \overset{\partial^0}{\longrightarrow} 
\mathcal{C}^1 \overset{\partial^1}{\longrightarrow}
\mathcal{C}^2 \overset{\partial^2}{\longrightarrow} \cdots
\overset{\partial^{m-2}}{\longrightarrow} \mathcal{C}^{m-1}
\overset{\partial^{m-1}}{\longrightarrow} \mathcal{C}^m 
\overset{\partial^m}{\longrightarrow} 0.
\end{equation}
Based on these definitions, the {\sl $k$-th cohomological space
$H^k(\frak g,V)$} is defined to be the following quotient space:
\[
H^k\big(\mathfrak{g},V\big) 
= 
\frac{\mathcal Z^k(\frak g,V)}{\mathcal B^k(\frak g,V)},
\]
where $\mathcal Z^k(\frak g,V):={\rm ker} \big( \partial^k \big)$ and
$\mathcal B^k(\frak g,V) := {\rm im} \big( \partial^{k-1} \big)$.

Within {\sc Maple}, there exists a package entitled
\textsf{LieAlgebraCohomology} which computes a somewhat different type
of Lie algebra cohomology, called {\sl relative cohomology}. In
particular, this package computes the {\sl De Rham cohomoloy}, quite
central in differential geometry. But still, there is no package or
command for computing the above-mentioned type of cohomological spaces
of Lie (super) algebras, although it has several applications to, {\em
e.g.}, the differential geometry of Cartan-Tanaka connections.

The article is divided in five sections. In Section~2, as already
said, some preliminaries about Gr\"obner bases are
reminded. Section~\ref{Computation} is devoted to the main results of
this paper. In Section~\ref{Description} we describe our algorithm to
compute the cohomological spaces of certain Lie algebras. Lastly, in
Section~\ref{Improved} we
show, with some examples, that computations naturally 
split up when the graduations are available. 

\section{Gr\"obner Bases and Elimination Ideals}
\label{Groebner} 

The theory of Gr\"obner bases is a key computational tool for studying
polynomial ideals. This theory was introduced and developed by
Buchberger, who devised its general scheme in the early 1960's
(\cite{Bruno2, Bruno3}). Nowadays, there exist several refined and
improved algorithms that are more efficient than the original one,
such as {\rm F}$_4$, {\rm F}$_5$, {\rm G$^2$V} and {\rm GVW}, and most
of them have been regularly implemented in computer algebra systems
like {\sc Maple}, {\sc Magma}, {\sc Mathematica}, {\sc Singular}, {\sc
Macaulay2}, {\sc Cocoa} and {\sc Sage}.

To provide a summarized description of the theory, 
let $\mathbb{ K} [x_1, \ldots,
x_n]$ be a polynomial ring in $n \geqslant 1$ variables
on some arbitrary commutative field
$\mathbb{K}$ of characteristic zero
and let $\mathcal{ I} = \langle f_1,\ldots ,f_k\rangle$
be any ideal of $\mathbb{ K} [x_1, \ldots, x_n]$ generated by a finite
number (noetherianity!) of polynomials $f_1,\ldots ,f_k\in \mathbb{ K}
[x_1, \ldots, x_n]$.

\begin{Definition}
A {\sl monomial ordering} on $\mathbb{ K}[ x_1, \dots, x_n]$ is a
relation $\prec$ on the set of monomials $x^\alpha = 
x_1^{ \alpha_1} \cdots x_n^{ \alpha_n}$ in $\mathbb{ K}
[x_1, \dots, x_n]$ which satisfies:

\begin{itemize}

\item $\prec$ is a total ordering;

\item $x^\alpha \prec x^\beta$ implies $x^\gamma x^\alpha \prec
x^\gamma x^\beta$ for every monomial $x^\gamma$, $\gamma \in 
\mathbb{ N}^n$;

\item $\prec$ is a well ordering.

\end{itemize}
\end{Definition}

For example, the usual {\sl lexicographical} ordering, 
here denoted $\prec_{\sf
lex}$, is a monomial ordering defined as follows (\cite{Becker,
Little}): if $\deg_i(m)$ denotes the degree in $x_i$ of a monomial
$m$, if $m'$ and $m''$ are two monomials, then $m' \prec_{\sf lex}
m''$ if and only if (by definition) the first nonzero entry of the
vector of $\mathbb{ Z}^n$:
\[
\big(\deg_1(m'')-\deg_1(m'),
\ldots, 
\deg_n(m'')-\deg_n(m')\big)
\]
is positive.

Let now $\prec$ be any monomial ordering on $\mathbb{ K} [x_1, \dots,
x_n]$. The {\sl leading monomial} of a polynomial $f \in \mathbb{ K}
[x_1, \ldots, x_n]$ is the greatest monomial\,\,---\,\,with respect to
$\prec$\,\,---\,\,which appears in $f$, and we denote it by ${\rm
LM}(f)$. Furthermore, the {\sl leading coefficient} of $f$, written by
${\rm LC}(f) \in \mathbb{ K}$, is the 
$\K$-coefficient of ${\rm LM}(f)$ in
$f$ and the {\sl leading term} of $f$ is the complete thing:
\[
{\rm LT}(f)
:= 
{\rm LC}(f)\cdot {\rm LM}(f). 
\]
The following theorem states a fundamental 
{\sl division algorithm} in $\mathbb{ K} [x_1, \ldots, x_n]$.

\begin{Theorem}
\label{division-Groebner}
{\rm (\cite{ Becker, Little})}
Given a fixed monomial ordering $\prec$ on $\mathbb{ K} [x_1, \dots,
x_n]$, for any {\em ordered} $k$-tuple $(f_1,\ldots,f_k)$ of polynomials
in $\mathbb{ K} [x_1, \ldots, x_n]$, every $f\in \mathbb{ K} [x_1,
\ldots, x_n]$ can be written as:
\[
f
=
a_1f_1+\cdots+a_kf_k+r,
\]
for some $a_i,r \in \mathbb{ K} [x_1, \ldots, x_n]$, with the main
property that either $r=0$ or $r$ is a linear combination of
monomials, {\em none of which} is divisible by any ${\rm LT}(f_j)$, $j
= 1, \dots, k$.
\end{Theorem}

Usually, one calls $r$ a {\em (one)} remainder of $f$ on division by
$(f_1, \dots, f_k)$, because most often, it is {\em not} unique, and
because in addition, it strongly depends on the ordering of the
$f_i$'s. This theorem, a higher-dimensional version of the standard
Euclidean division algorithm valid for the one-dimensional ring $\mathbb{
K} [ x_1]$, is the main effective cornerstone in the field of
Gr\"{o}bner bases; in fact, search for higher speed concentrates
mainly on improving the efficiency of division. Next, we define what
is a Gr\"obner basis for a polynomial ideal $\mathcal{ I} \subset
\mathbb{ K} [x_1, \ldots, x_n]$.

\begin{Definition}
A finite subset ${\tt G} = \{g_1, \ldots, g_l\} \subset \mathcal{ I}$ is
called a {\sl Gr\"obner basis} of $\mathcal{ I}$ with respect to some
fixed monomial ordering $\prec$ if the ideal generated by the leading
monomials of all elements of $\mathcal{ I}$ coincides with the
monomial ideal generated by the ${\rm LT} ( g_j)$, $j = 1, \dots, l$:
\[
\big<\,
{\rm LT}(f)\colon f\in\mathcal{I}\,\big> 
= 
\big<\,
{\rm LT}(g_1),\ldots,{\rm LT}(g_l)
\,\big>.
\]
\end{Definition}

Next, if ${\tt G} = \{g_1, \ldots, g_l\}$ is a Gr\"obner basis of an
ideal with respect to some monomial ordering $\prec$, one proves that
the remainder, on division by ${\tt G}$, of any $f\in \mathbb{ K} [
x_1, \dots, x_n]$ is {\em unique}, one calls this remainder the {\sl
normal form} of $f$ with respect to ${\tt G}$ and one denotes it by
${\rm NF}_{\tt G}(f)$, {\em cf.} again~\cite{ Becker, Little}. Also,
one proves that if ${\tt G}$ is a Gr\"{o}bner basis then ${\rm
NF}_{\tt G}(f)=0$ if and only if $f \in \left<\, {\tt G}
\,\right>$ belongs to the
ideal $\left<\, {\tt G} \,\right> =
\mathcal{ J}$. Then the fundamental theorem of the theory is that {\em
every} nonzero ideal $\mathcal{ I} \subset \mathbb{ K} [ x_1, \dots,
x_n]$ possesses at least one Gr\"obner basis, with (refinable)
algorithms which produces such a Gr\"obner
basis from any set of generators, by
taking so-called {\sl $S$-polynomials} between
any two distinct generators and by applying, 
inductively, the division Theorem~\ref{division-Groebner}. 
Furthermore, if ${\tt
G}$ is any Gr\"obner basis of $\mathcal{ I}$, it also generates
$\mathcal{ I}$, hopefully. However, Gr\"obner bases for an ideal are
not unique. Once a monomial order is chosen,
reduced Gr\"obner bases fully insure uniqueness.

\begin{Definition}
A {\sl reduced Gr\"obner basis} of an ideal $\mathcal{ I}$ is a
Gr\"obner basis ${\tt G} = \{ g_1, \dots, g_l\}$ of $\mathcal{ I}$ whose
polynomials $g_j$ are all monic such that, for any two distinct $g_{
j_1}, g_{ j_2} \in {\tt G}$, no monomial appearing in $g_{ j_2}$ is a
multiple of ${\rm LT}(g_{ j_1})$.
\end{Definition}

Then one establishes (\cite{Becker, Little}) that, given a fixed
monomial ordering $\prec$ on the ring $\mathbb{ K} [ x_1, \dots,
x_n]$, every ideal $\mathcal{ I} \subset \mathbb{ K} [ x_1, \dots,
x_n]$ possesses a {\em unique} reduced Gr\"obner basis.

\smallskip 
The concept of {\sl elimination ideal}, a natural application of
Gr\"obner bases, will be a very useful tool for us. Consider again
$\mathbb{K}[x_1,\ldots,x_n]$ and pick a (finite) subset of $m$, with
$1 \leqslant m \leqslant n-1$, variables among the $n$ variables
$\{x_1,\ldots,x_n\}$; possibly after a permutation, these
(sub)variables may of course be assumed to be just $x_1, \dots, x_m$.
Then, for any ideal $\mathcal{ I} \subset \mathbb{ K} [x_1, \dots,
x_m, x_{ m+1}, \dots, x_n]$, we call:
\[
\mathcal{I}\cap\mathbb{K}[x_1,\dots,x_m], 
\]
the {\sl elimination ideal} of $\mathcal{ I}$ with respect to the
(sub)variables:
\[
\big\{x_1,\dots,x_m\}
\subset
\big\{
x_1,\dots,x_m,x_{m+1},\dots,x_n
\big\}.
\]
The following proposition provides one with a way to compute elimination
ideals, using Gr\"obner bases, and, as a bonus, it also
yields at the same time a reduced Gr\"obner
basis for the elimination ideal.

\begin{Proposition}
\label{elimination-ideal}
{\rm (\cite{Becker, Little})}
Let $\prec$ be a monomial ordering on the ring $\mathbb{ K} [x_1,
\dots, x_m, x_{ m+1}, \dots, x_n]$ having the property that $x_j \prec
x_k$ for any $j = 1, \dots, m$ and any $k = m+1, \dots, n$, and let
${\tt G}$ be the reduced Gr\"obner basis of $\mathcal{ I}$ with
respect to $\prec$. Then ${\tt G} \cap \mathbb{K}[x_1, \dots, x_m]$ is
a reduced Gr\"obner basis for the elimination ideal $\mathcal{ I} \cap
\mathbb{ K}[x_1, \dots, x_m]$ with respect to $\prec$.
\end{Proposition}

Using this proposition, computers provide without pain\,\,---\,\,when
calculations succeed\,\,---\,\,elimination ideals, thanks to the
strength of implemented Gr\"obner bases. In particular, this gives a
simple way to solve systems of polynomial equations, even when they
have infinitely many solutions, and here presently, we shall have to
deal with solutions of equations that are {\em linear}, 
a case where calculations do most often succeed indeed.

\section{Computation of Cohomology Spaces}
\label{Computation}

Now, coming back to our goal, let $\frak g=\frak
g_{\overline{0}}\oplus\frak g_{\overline{1}}$ be an $m$-dimensional
Lie super algebra generated as a $\mathbb{ K}$-vector space
by $p$ even elements ${\sf e}_1,\ldots
{\sf e}_p$ and by $m-p$ odd elements ${\sf o}_{p+1},\ldots {\sf o}_m$,
and let $V$ be an $n$-dimensional $\frak g$-module generated by vectors
$v_1,\ldots,v_n$, as a
$\mathbb{ K}$-vector space too. It is natural
to divide any algorithm on the computation of Lie super algebra
cohomologies into three steps:

\begin{itemize}

\smallskip\item[$\bullet$] 
computation of the {\sl space of cocycles} ${\mathcal Z}^k(\frak g,V)$;

\smallskip\item[$\bullet$] 
computation of the {\sl space of coboundaries} 
${\mathcal B}^k(\frak g,V)$;

\smallskip\item[$\bullet$] 
computation of the cohomology space $H^k(\frak g,V) = 
{\mathcal Z}^k(\frak g,V) \big/ {\mathcal B}^k(\frak g,V)$.

\end{itemize}\smallskip

Sometimes, we shall abbreviate simply by ${\mathcal Z}^k$ the space
${\mathcal Z}^k(\frak g,V)$, and so on. Obviously, the most
substantial step of the algorithm is the third one, in which one has
to compute the quotient of the two spaces obtained, at the first and
second steps, by somewhat routine computations. Accordingly, we shall
divide this section into three steps in which we explain the
corresponding fraction of the algorithm.

\subsection{Computation of $\mathbf{ {\mathcal Z}^k (\frak g,V)}$}
At first, we have to determine a basis for the vector space $\mathcal
C^k(\frak g,V)$. For any $r = 0, \dots, k$, for any $1\leqslant i_1 <
\cdots < i_r \leqslant p$, for any $p+1 \leqslant j_{r+1} < \cdots <
j_k \leqslant m$ and for any $l = 1, \dots, n$, let us denote by:
\[
\Lambda^{(i_1,\ldots,i_r|j_{r+1},\ldots,j_k)}_l
\] 
the basic element (map)
of $\mathcal C^k (\frak g,V)$ whose value on $( {\sf
e}_{i_1}, \ldots, {\sf e}_{i_r}, {\sf o}_{ j_{r+1}}, \ldots, {\sf
o}_{j_k})$ is exactly $1\cdot v_l$, which acts super-symmetrically and
which is zero elsewhere. One verifies that the set of these $n\, {m
\choose k}$ maps constitutes a basis over $\mathbb{ K}$ for the vector
space $\mathcal C^k (\frak g,V)$, hence a general $k$-cochain $\Phi$
naturally decomposes as a linear combination:
\[
\Phi
=
\sum_{r=0}^k\,
\sum_{1\leqslant i_1 < \cdots < i_r \leqslant p}\,
\sum_{ p+1 \leqslant j_{r+1} < \cdots < j_k\leqslant m}\,
\sum_{l=1}^n\,
\phi_{(i_1,\ldots,i_r|j_{r+1},\ldots,j_k)}^l\ 
\Lambda^{(i_1,\ldots,i_r|j_{r+1},\ldots,j_k)}_l,
\]
where the $\phi_{ (i_1, \ldots, i_r | j_{r+1}, \ldots, j_k)}^l \in
\mathbb{ K}$ are arbitrary scalars in the ground field. For more brevity
and without much abuse of notation, let us denote $\phi_{ (i|j)_{
r,k}}^l$, $\Lambda^{ (i|j)_{r,k}}_l$ and $({\sf e}_i,{\sf o}_j)_{r,k}$
instead of $\phi_{ (i_1, \ldots, i_r|j_{r+1}, \ldots, j_k)}^l$,
$\Lambda^{ (i_1, \ldots, i_r|j_{r+1}, \ldots, j_k)}_l$ and $({\sf
e}_{i_1}, \ldots, {\sf e}_{i_r}, {\sf o}_{ j_{r+1}}, \ldots,{\sf o}_{
j_k})$, respectively. Thus, with these abbreviated notations, the above
expansion of a general $k$-cochain reads:
\begin{equation}
\label{Phi}
\Phi
=
\sum_r\,
\sum_{i_1<\cdots<i_r}\,\sum_{j_{r+1}<\cdots<j_k}\,\sum_l\,
\phi_{(i|j)_{r,k}}^l\,\,
\Lambda^{(i|j)_{r,k}}_l.
\end{equation}
In the important (special) case of standard Lie algebras
$\mathfrak{ g} \subset \mathfrak{ h} = V$ represented
by means of bases: 
\[
\mathfrak{g}
=
\K\,{\sf e}_1
\oplus\cdots\oplus
\K\,{\sf e}_m
\ \ \ \ \ 
\text{\rm and}
\ \ \ \ \
\mathfrak{h}
=
\K\,{\sf f}_1
\oplus\cdots\oplus
\K\,{\sf f}_n,
\]
odd elements
are plainly absent, whence the expression of a general
$k$-cochain reduces to:
\[
\Phi
=
\sum_{1\leqslant i_1<\cdots<i_k\leqslant m}\,
\sum_{l=1}^n\,\,
\phi_{i_1,\ldots,i_k}^l\,\, 
\Lambda^{i_1,\ldots,i_k}_l,
\]
where the basic $k$-cochains $\Lambda_l^{ 
i_1, \dots, i_k}$ 
also write as follows 
in terms of the dual ${\sf e}_i^*$:
\[
\Lambda_l^{i_1,\dots,i_k}
=
{\sf e}_{i_1}^*\wedge\cdots\wedge{\sf e}_{i_k}^*
\otimes
{\sf f}_l.
\]

Now, in order to compute the cocycle subspace ${\mathcal Z}^k \subset
\mathcal{ C}^k$, one proceeds by applying the fundamental
formula~\thetag{ \ref{cohomology}} to know what value $\partial^k
\Phi$ has on each $(k+1)$-tuple $({\sf e}_i,{\sf o}_j)_{ s, k+1}$, for
all $s = 0, \dots, k+1$, for all $1 \leqslant i_1< \cdots <i_s
\leqslant p$, for all $p+1 \leqslant j_{s+1}< \cdots < j_{
k+1}\leqslant m$, and afterwards, by just equating to zero each such
expression $( \partial^k \Phi) \big(( {\sf e}_i, {\sf o}_j)_{ s,k+1}
\big)$, a task which is of course left to a computer. With more
precisions, because each such $(\partial^k \Phi) \big(( {\sf e}_i,
{\sf o}_j)_{ s,k+1} \big)$ belongs to the $n$-dimensional $\mathbb{
K}$-vector space $V$, one in fact gets $n$ scalar equations in this
way. After all, this gives in sum 
exactly $n\, {m \choose k+1}$ homogeneous equations
that are all linear with respect to the $n\, {m \choose k}$ unknown
coefficients $\phi_{(i|j)_{r,k}}^l$.  Then by computer-solving the
obtained linear system which we shall denote by: 
\[
{\sf Syst}_\phi\big(\mathcal{ Z}^k\big), 
\]
one completely identifies those coefficients $\phi_{ (i|j)_{ r,k}}^l$
which make up cocycles $\Phi = \sum\, \phi_{ (i|j)_{ r,k}}^l\,
\Lambda_l^{ (i | j)_{ r, k}}$ which belong to ${\mathcal Z}^k$. The first
step ends so.

\subsection{Computation of $\mathbf{{\mathcal B}^k(\frak g,V)}$}
This second step is rather similar to the first one, though
less direct, for
it requires the use of elimination ideals
(Proposition~\ref{elimination-ideal}). Indeed using
once more the general representation~\thetag{ \ref{Phi}}
with $k$ replaced by $k-1$, a general $(k-1)$-cochain 
writes quite similarly under the form:
\begin{equation}
\label{Psi}
\Psi
=
\sum_{r=0}^{k-1}\,
\sum_{1\leqslant i_1<\cdots<i_r\leqslant p}\,
\sum_{p+1\leqslant j_{r+1}<\cdots<j_{k-1}\leqslant m}\,
\sum_{l=1}^n\,
\psi_{(i|j)_{r,k-1}}^l\,\,
\Lambda^{(i|j)_{r,k-1}}_l,
\end{equation}
where the $\psi_{(i|j)_{r,k-1}}^l \in \mathbb{ K}$ are arbitrary
scalars in the ground field. By definition, the elements of $\mathcal
B^k$, namely the coboundaries, are $k$-cochains of the form
$\partial^{ k-1} \Psi$, for such a $\Psi$.  With more precision,
$\mathcal B^k$ is the space of $k$-cochains $\Phi$ as in~\thetag{
\ref{Phi}} that are of the form $\Phi = \partial^{ k-1} \Psi$, for
some $(k-1)$-cochains $\Psi$ as in~\thetag{ \ref{Psi}}.  Consequently,
applying once again the fundamental formula~\thetag{
\ref{cohomology}}, we have to compute the value of $\partial^{ k-1}
\Psi$ on each of the $k$-tuples $({\sf e}_i,{\sf o}_j)_{r,k}$
belonging to $\frak g^k$ and then to equate them to the value of
$\Phi$ on these $k$-tuples, where we recall that:
\[
\Phi
\big(({\sf e}_i,{\sf o}_j)_{r,k}\big)
=
\Phi({\sf e}_{i_1},\ldots,{\sf e}_{i_r},{\sf
o}_{j_{r+1}},\ldots,{\sf o}_{j_k})
=
\sum_{l=1}^n\,
\phi_{(i_1,\ldots,i_r|j_{r+1},\ldots,j_k)}^l \
v_l.
\]
But looking at~\thetag{ \ref{cohomology}}, and without performing
explicit computations (left to a computer in specific examples), one
easily convinces oneself that there are certain {\em linear} forms
${\sf L}_{ i, j, r, k}$ in the coefficients
$\psi_{(i'|j')_{r',k-1}}^{l'}$ of $\Psi$ such that:
\[
(\partial^{k-1}\Psi)
\big(({\sf e}_i,{\sf o}_j)_{r,k}\big)
=
\sum_{l=1}^n\,
{\sf L}_{i,j,r,k}
\big(
\big\{\psi_{(i'|j')_{r',k-1}}^{l'}\big\}
\big)\,
v_l.
\]
Hence for any $i$, $j$, $r$, $k$, by equating the coefficients of the
$v_l$, $l = 1, \dots, n$, in both sides of the equalities:
\[
\partial^{k-1}\Psi
\big(({\sf e}_i,{\sf o}_j)_{r,k}\big)
=
\Phi\big(({\sf e}_i,{\sf o}_j)_{r,k}\big),
\]
it therefore follows that a $k$-cochain $\Phi = \partial^{ k-1} \Psi$
is a $k$-coboundary if and only if all its coefficients $\phi_{ (i |
j)_{ r, k}}^l$ are of the form:
\[
\phi_{(i|j)_{r,k}}^l
=
{\sf L}_{i,j,r,k}
\Big(
\big\{\psi_{(i'|j')_{r',k-1}}^{l'}\big\}
\Big), 
\]
for {\em some} $(k-1)$-cochain $\Psi$ having coefficients
$\psi_{(i'|j')_{r',k-1}}^{l'}$. The task of writing explicitly the
right-hand sides being left to a computer, we obtain in this way $n\,
{m\choose k}$ linear equations. Lastly, we can use Gr\"obner bases to
{\em eliminate} all the variables $\psi_{(i'|j')_{r',k-1}}^{ l'}$ in
these linear equations ({\em cf.}
Proposition~\ref{elimination-ideal}), which provides at the end a
collection of linear equations (automatically organized as a reduced
Gr\"obner basis) involving only the variables
$\phi_{(i|j)_{r,k}}^l$. If we denote this new system by:
\[
{\sf Syst}_\phi\big(\mathcal{B}^k\big), 
\]
the fact that one always has $\mathcal{ B}^k \subset \mathcal{ Z}^k$
entails that any solution of ${\sf Syst}_\phi \big( \mathcal{ B}^k
\big)$ is necessarily a solution of ${\sf Syst}_\phi \big( \mathcal{
Z}^k \big)$. However as usual in linear algebra, this does not mean
that the (finite) collection of equations for ${\sf Syst}_\phi \big(
\mathcal{ Z}^k \big)$ is {\em included}, as a set, in the (finite)
collection of equations for ${\sf Syst}_\phi \big( \mathcal{ B}^k
\big)$: one in general needs to make linear combinations until this
becomes true.

\subsection{Computation of $\mathbf{H^k(\frak g,V)}$}
Now we are ready to start the third, main step, namely the computation
of the $k$-th cohomological space $H^k = \mathcal{ Z}^k \big/
\mathcal{ B}^k$.  (Of course, any technique which decreases the
complexity of this last step simultaneously increases the speediness
of computations.)  The two systems ${\sf Syst}_\phi ( \mathcal{ Z}^k)$
and ${\sf Syst}_\phi ( \mathcal{ B}^k)$ of linear equations in the
unknown variables $\phi_{ (i|j)_{ r,k}}^l$ identify exactly all the
elements of $\mathcal Z^k$ and $\mathcal B^k$,
respectively. Therefore, every nonzero element of the quotient
$\mathbb{ K}$-vector space:
\[
H^k
=
\mathcal{Z}^k
\big/
\mathcal{B}^k
=
\mathcal{Z}^k\,\,
{\rm mod}\,\mathcal{B}^k
\]
is of the form: 
\[
\Phi 
+ 
\mathcal{B}^k, 
\]
where the coefficients $\phi^{(i|j)_{r,k}}_l$ of the $k$-cochain $\Phi
= \sum\, \phi_{ (i | j)_{ r, k}}^l\, \Lambda_l^{ (i | j)_{ r, k}}$
satisfy all the equations in ${\sf Syst}_\phi ( \mathcal{ Z}^k)$ and
do not satisfy at least one of the equations in ${\sf Syst}_\phi (
\mathcal{ B}^k)$.

\subsection{Finding a basis for a quotient $\K$-vector space}
Temporarily, let us set aside our cohomological objective and let us
present some results in the theory of Gr\"{o}bner basis that are
useful to the purpose of finding representatives of the quotient
$V / W = V \, {\rm mod}\, W$ of
any two $\mathbb{ K}$-vector subspaces $W \subset V \subset E$ sitting
inside a certain (large) ambient $\K$-vector space $E$. 

In a first moment, 
given a vector subspace $F \subset E$ of some $\K$-vector space $E$
which is represented as the zero-set of some linear forms\,\,---\,\,as
for instance $\mathcal{ Z}^k \subset \mathcal{ C}^k$ which is
represented by ${\sf Syst}_\phi \big( \mathcal{ Z}^k \big)$\,\,---, by
allowing fully the use of Gr\"obner bases, we want to
find an explicit set of vectors ${\sf f}_1, \dots, {\sf f}_{\dim F}
\in E$ which make up a basis for $F$. Then in a second moment and still
employing Gr\"obner bases, given instead two $\K$-vector subspaces $W
\subset V \subset E$ of dimensions $p := \dim_\K V$ and $q := \dim_\K
W$ which are both represented as zero-sets of some linear
forms\,\,---\,\,as for instance $\mathcal{ B}^k \subset \mathcal{ Z}^k
\subset \mathcal{ C}^k$ which are represented by ${\sf Syst}_\phi
\big( \mathcal{ B}^k \big)$ and by ${\sf Syst}_\phi \big( \mathcal{
Z}^k \big)$\,\,---, we will show how to find explicitly $p - q$
linearly independent vectors ${\sf v}_1, \dots, {\sf v}_{ p - q} \in
V$ such that:
\[
{\sf v}_1
+W,
\,\,\dots,\,\,
{\sf v}_{p-q}+W
\]
make up a basis for the quotient vector space $V / W = V\, {\rm mod}
\, W$.

Thus, let $E$ be a $\K$-vector space of dimension $n \geqslant 1$, let
$\{ {\sf e}_1, \dots, {\sf e}_n\}$ be a basis of $E$ and let $(x_1,
\dots, x_n) \in \K^n$ be the associated coordinates in terms of which
any vector ${\sf e} \in E \simeq \K^n$ represents uniquely as:
\[
{\sf e}
=
x_1\,{\sf e}_1
+\cdots+
x_n\,{\sf e}_n.
\]
By convention, the variable names $x_i$ will be reserved to write down
Cartesian equations of vector subspaces, and we will
also need some other auxiliary variables $(y_1, \dots, y_n)$.

To begin with, consider the circumstance where a
given vector subspace
$F \subset E \simeq \K^n$ is represented 
as generated by $\mu$ vectors
${\sf f}_1, \dots, {\sf f}_\mu \in F$ that are not necessarily
linearly independent. Each such vector decomposes
according to the basis:
\[
{\sf f}_1
=
f_{11}\,{\sf e}_1
+\cdots+
f_{1n}\,{\sf e}_n,
\,\,\dots\dots,\,\,
{\sf f}_\mu
=
f_{\mu1}\,{\sf e}_1
+\cdots+
f_{\mu n}\,{\sf e}_n,
\]
for some scalars $f_{ \lambda i} \in \K$, and using the
auxiliary variables $(y_1, \dots, y_n)$, 
we associate to them the following $\mu$ linear forms:
\[
f_1(y)
:=
f_{11}\,y_1
+\cdots+
f_{1n}\,y_n,
\,\,\dots\dots,\,\,
f_\mu(y)
:=
f_{\mu 1}\,y_1
+\cdots+
f_{\mu n}\,y_n,
\]
which we simply view as (degree $1$) {\em polynomials} belonging to
$\K [ y_1, \dots, y_n]$. 
The proofs of the three statements below,
including the following preliminary proposition, will be postponed to
the end of the present section.

\begin{Proposition}
\label{f-mu-g-m}
Fix a lexicographic ordering $\prec$ on monomials of the ring $\K[
y_1, \dots, y_n]$. With $F = {\rm Vect}_\K ( {\sf f}_1, \dots, {\sf
f}_\mu)$ as above, and with the associated linear forms $f_1 ( y),
\dots, f_\mu ( y)$, if ${\tt G} := \{ g_1 ( y), \dots, g_m(y) \}$ is
{\em the} reduced Gr\"obner basis of the ideal:
\[
\big<\,
f_1(y),\dots,f_\mu(y)
\,\big>
\]
in $\K[ y_1, \dots, y_n]$ with respect to $\prec$, then:

\begin{itemize}

\smallskip\item[{\bf (i)}]
$\dim_\K F = m =$ precisely the cardinal of ${\tt G}${\em ;}

\smallskip\item[{\bf (ii)}]
all $g_j (y)$, $j = 1, \dots, m$, are linear forms, namely:
\[
g_j(y)
=
g_{j1}\,y_1
+\cdots+
g_{jn}\,y_n
\]
for some scalars $g_{ ji} \in \K$, and furthermore, the
$m$ vectors:
\[
{\sf g}_1
:=
g_{j1}\,{\sf e}_1
+\cdots+
g_{jn}\,{\sf e}_n,
\,\,\dots\dots,\,\,
{\sf g}_m
:=
g_{m1}\,{\sf e}_1
+\cdots+
g_{mn}\,{\sf e}_n
\]
constitute a basis for $F$ as a vector space;

\smallskip\item[{\bf (iii)}]
an arbitrary vector ${\sf h} = h_1\, {\sf e}_1 + \cdots + h_n\, {\sf
e}_n \in E$, with coordinates $h_i \in \K$, belongs to $F$ {\em if and
only if} the normal form of the
associated $h(y) := h_1\, y_1 + \cdots + h_n\, y_n$
with respect to the reduced Gr\"obner basis ${\tt G}$ is zero:
\[
0
=
{\rm NF}_{\tt G}(h).
\]
\end{itemize}

\end{Proposition}

However, as we said, the $\K$-vector subspace $F \subset E$ we want to
consider for applications to (super) Lie algebra cohomologies, namely
$\mathcal{ Z}^k \subset \mathcal{ C}^k$ (or also $\mathcal{ B}^k
\subset \mathcal{ C}^k$) should be thought of as being represented as
the zero-set of some (Cartesian) linear equations. The appropriate
statement will better be brought to light by means of a simple
illustration.

\begin{Example}
Consider the system of three (Cartesian) linear equations:
\[
\left\{ 
\aligned
f_1(x)
&
:= 
x_1-x_4+x_5=0,
\\
f_2(x)
&
:= 
2\,x_1+x_2+x_4=0,
\\
f_3(x)
&
:=
-x_3+2\,x_4+x_5=0,
\endaligned\right.
\]
in the vector space $E = \K^5$
with coordinates $(x_1, x_2, x_3, x_4, x_5)$ which
represents a certain vector subspace $F \subset E$. Transforming
(either by hand or with a computer) the ideal $\big<\, f_1 ( x), \,
f_2 ( x), \, f_3 ( x)\, \big>$ to the reduced Gr\"obner basis with
respect to the lexicographic ordering $x_5 \prec x_4 \prec x_3 \prec
x_2 \prec x_1$, one gets that $F \subset E$ is equivalently defined as
the set of all $(x_1,x_2,x_3,x_4,x_5) \in \K$ satisfying: $0 = g_1(x)
= g_2(x) = g_3 (x)$, where:
\[
g_1(x)
:=
x_1-x_4+x_5,\ \ \ \ \ \ \
g_2(x)
:=
x_2+3\,x_4-2\,x_5,\ \ \ \ \ \ \
g_3(x)
:=
x_3-2\,x_4-x_5,
\]
and where ${\tt G} := \{ g_1 ( x), g_2 ( x), g_3 ( x)\}$ is the
reduced Gr\"obner basis in question.  Thus, $x_4$ and $x_5$, are {\em
horizontal parameters} for $F$, $x_1, x_2, x_3$ are functions of
$(x_4, x_5)$, and $F$ is a {\em graphed}\, $5-3=2$-dimensional
subspace of the $5$-dimensional vector space $E = \K^5$.

Next, choosing firstly $(x_4 , x_5) = (1, 0)$ and secondly $(x_4, x_5)
= (0, 1)$, one sees that $F$ is generated by the two column vectors
$(1,-3,2,1,0)^{\sf t}$ and $(-1,2,1,0,1)^{\sf t}$. To these two
vectors, one then associates the following set of two linear forms:
\[
\Big\{y_1-3y_2+2y_3+y_4,\,\,-y_1+2y_2+y_3+y_5\Big\},
\]
in some five auxiliary variables $y_1, y_2, y_3, y_4, y_5 \in \K$.  On
the other hand, granted that computing a normal form with respect to
${\tt G}$ just means replacing $x_1$ by $x_4 - x_5$, $x_2$ by $-3\,
x_4 + 2\, x_5$ and $x_3$ by $2\, x_4 + x_5$, and considering
the auxiliary bilinear form $\sum_{ i=1}^5 \, x_i \, y_i$, 
we see that:
\[
{\rm NF}_{\tt G}
\Big(
{\textstyle{\sum_{i=1}^5}}\,
x_i\,y_i
\Big)
=
(x_4-x_5)\,y_1+(-3\,x_4+2\,x_5)\,y_2
+
(2\,x_4+x_5)\,y_3+x_4\,y_4+x_5\,y_5.
\]
Reorganizing, we easily find the coefficients of the parameters $x_4$
and $x_5$ in this expression:
\[
\aligned
&
\boxed{x_4}\,\colon \ \ \ \ \ \ y_1-3\,y_2+2\,y_3+y_4
\\
&
\boxed{x_5}\,\colon \ \ \ \ \ \ -\,y_1+2\,y_2+y_3+y_5,
\endaligned
\]
and interestingly enough, these two coefficients coincide with the
above two linear forms in the auxiliary variables $y_i$. This is a
quite general fact, whose proof is also postponed to the end of the
present section.
\end{Example}

\begin{Proposition}
\label{extract-basis-from-Cartesian}
Let $F \subset E \simeq \K^n$ be a $\K$-vector subspace which is
represented by means of Cartesian linear equations:
\[
F
=
\big\{
\text{\rm vectors}\,\,
x_1\,{\sf e}_1
+\cdots+
x_n\,{\sf e}_n\,\,
\text{\rm s.t.}\,\,
0
=
f_1(x)
=\cdots=
f_\mu(x)
\big\},
\]
for a certain collection of $\mu \geqslant 1$ linear forms
$f_\lambda ( x)$. 
Let ${\tt G}$ be the reduced Gr\"obner basis of the ideal
$\big< \, f_1 ( x), \dots, f_\mu ( x)\, \big>$ with respect
to some fixed {\em lexicographic} ordering. Given $n$ new auxiliary
indeterminates $y_1, \dots, y_n$, let:
\[
h_y(x)
:=
{\rm NF}_{\tt G}
\big(
x_1y_1+\cdots+x_ny_n
\big)
\in
\K[x_1,\dots,x_n]
\]
be the normal form, with respect to ${\tt G}$, of
the bilinear form $\sum_{i=1}^n\, x_i\, y_i$. Then the following
four assertions hold true:

\begin{itemize}

\smallskip\item[{\bf (i)}]
$h_y(x)$ is linear in $(x_1, \dots, x_n)$;

\smallskip\item[{\bf (ii)}]
$h_y(x)$ involves exactly $\dim F =: m$ variables $x_i${\em :}
\[
h_y(x)
=
x_{i_1}\,h_1(y)
+\cdots+
x_{i_m}\,h_m(y),
\]
for some $1 \leqslant i_1 < \cdots < i_m \leqslant n${\em ;}

\smallskip\item[{\bf (iii)}]
all the appearing coefficients $h_j ( y)$ of $h_y (x)$
are linear forms in the
variables $(y_1, \dots, y_n)$; 

\smallskip\item[{\bf (iv)}]
if one expands them:
\[
h_j(y)
=
h_{j1}\,y_1
+\cdots+
h_{jn}\,y_n
\ \ \ \ \ \ \ \ \ \ \ \ \ 
{\scriptstyle{(j\,=\,1\,\cdots\,m)}}
\]
in terms of some scalars $h_{ ji} \in \K$, then the $m$ 
associated vectors:
\[
{\sf h}_1
:=
h_{11}\,{\sf e}_1
+\cdots+
h_{1n}\,{\sf e}_n,
\,\,\dots\dots,\,\,
{\sf h}_m
:=
h_{m1}\,{\sf e}_1
+\cdots+
h_{mn}\,{\sf e}_n
\]
make up a basis for $F$. 

\end{itemize}

\end{Proposition}

The last data ${\sf h}_1, \dots, {\sf h}_m$ are exactly what we
wanted: an explicit basis for the $\K$-vector subspace 
$F \subset E$ which was represented by linear
equations. 

\smallskip

We can now come back to our initial goal. Let $E \simeq \K^n$ be an
ambient
$n$-dimensional $\K$-vector space as above, fix coordinates $(x_1,
\dots, x_n)$ on $E$ and fix some {\em lexicographic} ordering on
monomials of $\K [ x_1, \dots, x_n]$.  Let $W \subset E$ and $V
\subset E$ be two $\K$-vector subspaces which are both represented by
means of Cartesian linear equations:
\[
\aligned
W
=
&
\big\{
\text{\rm vectors}\,\,
x_1\,{\sf e}_1
+\cdots+
x_n\,{\sf e}_n\,\,
\text{\rm s.t.}\,\,
0
=
g_1(x)
=\cdots=
g_\nu(x)
\big\},
\\
V
=
&
\big\{
\text{\rm vectors}\,\,
x_1\,{\sf e}_1
+\cdots+
x_n\,{\sf e}_n\,\,
\text{\rm s.t.}\,\,
0
=
f_1(x)
=\cdots=
f_\mu(x)
\big\},
\endaligned
\]
for certain two collections of linear forms $g_1 ( x), \dots, g_\nu
(x)$ and $f_1 ( x), \dots, f_\mu (x)$, with the further assumption
that $W \subset V$. For our cohomological objective, the initial data
are precisely presented under such form: $\mathcal{ B}^k \subset
\mathcal{ C}^k$ and $\mathcal{ Z}^k \subset \mathcal{ C}^k$ are the
zero-sets of ${\sf Syst}_\phi \big( \mathcal{ B}^k \big)$ and of ${\sf
Syst}_\phi \big( \mathcal{ B}^k \big)$, respectively, with $\mathcal{
B}^k \subset \mathcal{ Z}^k$, of course. It goes without
saying that
Proposition~\ref{extract-basis-from-Cartesian} provides two explicit
bases for $W$ and $V$, namely:
\[
\aligned
W
=
{\rm Span}_\K
\big(
{\sf w}_1,\dots,{\sf w}_q
\big)
\ \ \ \ \ \ \ \ \ \ 
\text{\rm and}
\ \ \ \ \ \ \ \ \ \ 
V
=
{\rm Span}_\K
\big(
{\sf v}_1,\dots,{\sf v}_p
\big),
\endaligned
\]
where $q := \dim_\K W$ and $p := \dim_\K V$. The following theorem
then realizes the goal of finding a basis for $V / W = 
V\, {\rm mod}\, W$ as a
$\K$-vector space.

\begin{Theorem}
\label{quotient-V-W}
Let $E$ be an $n$-dimensional $\K$-vector space equipped with a basis
$\{ {\sf e}_1, \dots, {\sf e}_n \}$, let $V \subset E$ and $W \subset
E$ be two $\K$-vector subspaces having dimensions $p := \dim_\K V$ and
$q := \dim_\K W$ that are both represented:
\[
\aligned
V
=
{\rm Span}_\K
\big(
{\sf v}_1,\dots,{\sf v}_p
\big)
\ \ \ \ \ \ \ \ \ \ 
\text{\rm and}
\ \ \ \ \ \ \ \ \ \ 
W
=
{\rm Span}_\K
\big(
{\sf w}_1,\dots,{\sf w}_q
\big),
\endaligned
\]
as the span of some basis vectors: 
\[
\aligned
{\sf v}_i
&
= 
v_{i1}\,{\sf e}_1
+\cdots+
v_{in}\,{\sf e}_n
\ \ \ \ \ \ \ \ \ 
\text{\rm and}
\ \ \ \ \ \ \ \ \ 
{\sf w}_j
=
w_{j1}\,{\sf e}_1
+\cdots+
w_{jn}\,{\sf e}_n
\\
&
\ \ \ \ \ \ \ \ \ \ \ \ \ 
{\scriptstyle{(i\,=\,1\,\cdots\,p)}}
\ \ \ \ \ \ \ \ \ \ \ \ \ \ \ \ \ \ \ \ \ \ \ \ \ \ 
\ \ \ \ \ \ \ \ \ \ \ \ \ \ \ \ \ \ \ \ \ \ \ \ \ \
{\scriptstyle{(j\,=\,1\,\cdots\,q)}}
\endaligned
\]
which are explicitly given in terms of their coordinates
$v_{ ik} \in \K$ and $w_{ jk}\in \K$. Suppose
that $W \subset V$, whence
$q \leqslant p$, and associate to these two bases the following
two collections of linear forms:
\[
\aligned
f_i(y)
:=
v_{i1}\,y_1
+\cdots+
v_{in}\,y_n
\ \ \ \ \ \ \ \ \ 
\text{\rm and}
\ \ \ \ \ \ \ \ \ 
g_j(y)
:=
w_{j1}\,y_1
+\cdots+
w_{jn}\,y_n
\endaligned
\]
in some auxiliary $\K [ y_1, \dots, y_n]$. 
Lastly, let: 
\[
{\tt B}_V
:=
\big\{
\overline{f}_1(y),\dots,\overline{f}_p(y)
\big\}
\ \ \ \ \ \ \
\text{\rm and} 
\ \ \ \ \ \ \
{\tt B}_W
:=
\big\{
\overline{g}_1(y),\dots,\overline{g}_q(y)
\big\}
\]
be the two reduced Gr\"obner bases of the two ideals $\big< \, f_1 ( y),
\dots, f_p ( y)\, \big>$ and $\big< \, g_1 ( y), \dots, g_q ( y)\,
\big>$ with respect to some fixed {\em lexicographic} ordering $\prec$
on
the monomials of $\K [ y_1, \dots, y_n]$. 
Then the reduced Gr\"obner
basis ${\tt B}_{V/W}$ of the ideal:
\[
\Big<\,
{\rm NF}_{{\tt B}_W}\big(\,\overline{f}\,\big)\,
\colon\,
\overline{f}\in{\tt B}_V
\,\Big>
\] 
generated by the normal forms with 
respect to ${\tt B}_W$ of all elements of ${\tt B}_V$,
is of cardinal equal to $p - q = \dim V - \dim W$, and furthermore,
if:
\[
\overline{h}_l(y)
=
h_{l1}\,y_1
+\cdots+
h_{ln}\,y_n
\ \ \ \ \ \ \ \ \ \ \ \ \ 
{\scriptstyle{(l\,=\,1\,\cdots\,p\,-\,q)}}
\]
are its elements, the $p - q$ associated vectors:
\[
{\sf h}_l
:=
h_{l1}\,{\sf e}_1
+\cdots+
h_{ln}\,{\sf e}_n
\,\,{\rm mod}\,W
\ \ \ \ \ \ \ \ \ \ \ \ \ 
{\scriptstyle{(l\,=\,1\,\cdots\,p\,-\,q)}}
\]
belong to $V$ and, when considered ${\rm mod}\, W$, 
make up a basis for $V / W = V\, {\rm mod}\, W$. 
\end{Theorem}

Computer tests ({\em cf.} examples below) show that, compared with
standard linear algebra methods, the use of Gr\"obner bases improves
speed and efficiency, especially because the computations underlying
Proposition~\ref{extract-basis-from-Cartesian} and
Theorem~\ref{quotient-V-W} can be achieved within a polynomial ring,
without the need of several transformations between polynomials and
vectors; indeed, from the two collections of Cartesian linear
equations ${\sf Syst}_\phi \big( \mathcal{ Z}^k \big)$ and ${\sf
Syst}_\phi \big( \mathcal{ B}^k \big)$,
Proposition~\ref{extract-basis-from-Cartesian} extracts two
collections of polynomials in some auxiliary variables $\upsilon_{ (i
\vert j)_{ r, k}}^l$ to which one can directly apply
Theorem~\ref{quotient-V-W} in order to find a basis for the sought
cohomology space $H^k = \mathcal{ Z}^k \big/ \mathcal{ B}^k$, {\em
see} also the description of the algorithm in the next section.

\proof[Proof of Proposition~\ref{f-mu-g-m}] We begin by making a
preliminary observation. According to the process of producing any
Gr\"obner basis, each element $g_j (y)$ of ${\tt G}$ is obtained by
subjecting all pairs $\{ f_{\lambda_1} (y), f_{\lambda_2}(y)\}$ to an
$S$-polynomial elimination of leading terms, by performing division
(Theorem~\ref{division-Groebner}) and by repeating the process until
stabilization, whence one easily convinces oneself that {\em only
linear forms}, namely degree one polynomials having no constant term,
can come up at each stage. At the end, every $g_j ( y)$ is therefore a
linear form. Of course, the ideal is the same:
\[
\left<\,f_1(y),\dots,f_\mu(y)\,\right>
=
\left<\,g_1(y),\dots,g_m(y)\,\right>.
\]
Thus, because all considered polynomials are linear forms, there
necessarily exist some scalars $c_{j \lambda} \in \mathbb K$ such that
$g_j (y) = \sum_{\lambda = 1}^\mu\, c_{ j \lambda} \, f_\lambda (y)$
for all $j = 1, \dots, m$, and in the other direction also, there
necessarily exist some scalars $d_{\lambda j} \in \mathbb K$ such that
$f_\lambda (y) = \sum_{j = 1}^m\, d_{\lambda j} \, g_j (y)$ for all
$\lambda = 1, \dots, \mu$.  It follows that the vector subspace
$F_{\tt G}$ associated to the $g_j$ by 
{\bf (ii)} is contained in the original
vector subspace $F \subset E$ to which the $f_\lambda (y)$ were
associated, and also in the other direction that $F \subset F_{\tt G}$.
Consequently, we have $F = F_{\tt G}$.

To finish with {\bf (i)} and {\bf (ii)}, it remains to prove the
linear independency of the vectors ${\sf g}_1, \dots, {\sf g}_m$
associated to $g_1 (y), \dots, g_m (y)$. Suppose by contradiction that
$0 = c_1 \, {\sf g}_1 + \cdots + c_m\, {\sf g}_m$ for some
$c_i\in\mathbb K$ that are not all zero.  It immediately follows that
$c_1 \, g_1 (y) + \cdots + c_m\, g_m (y) \equiv 0$. Consequently there
exist at least two different integers $j_1\neq j_2$ such that ${\rm
LM} (g_{ j_1}) = {\rm LM}(g_{ j_2})$, contrarily to the assumption
that the chosen ${\tt G}$ was a {\em reduced} Gr\"obner basis. In sum:
\[
m
=
{\rm Card}\,{\tt G} 
= 
\dim_\K F.
\]

Lastly, we check {\bf (iii)}. Of course, a vector ${\sf h}$ belongs to
$F = F_{\tt G}$ 
if and only if thre exist scalars $c_i \in \K$ such that ${\sf h}
= c_1\, {\sf g}_1 + \cdots + c_m\, {\sf g}_m$. Equivalently, the
associated polynomial (linear form) $h(y) = c_1 \, g_1 (y) + \cdots +
c_m \, g_m (y)$ belongs to the ideal generated by the Gr\"obner basis
${\tt G}$. But this is so if and only if the normal form ${\rm
NF}_{\tt G}(h)$ of $h(y)$ with respect to ${\tt G}$ is zero.
\endproof

\proof[Proof of Proposition~\ref{extract-basis-from-Cartesian}] We
already saw, in the beginning of the proof of the preceding
proposition, that all elements of ${\tt G}$ are linear forms and that
any division by ${\tt G}$ preserves linearity in $\K [x_1, \dots,
x_n]$. Since $\sum_{ i=1}^n\, x_i\, y_i$ is linear in the $x_i$, its
normal form $h_y (x)$ with respect to ${\tt G}$ is also linear, which
is {\bf (i)}.

Next, let $\underline{ m}$ denote the cardinal of the Gr\"obner basis
${\tt G}$ and denote its elements by $g_1 (x), \dots, g_{ \underline{
m}} (x)$.  Since ${\tt G}$ is reduced, for all $l = 1, \dots,
\underline{ m}$, the leading terms of $g_l (x)$ are monic, of degree
one of course, and distinct, say:
\[
x_{i_1}
=
{\rm LT}(g_1),
\dots,
x_{i_{\underline{m}}}
=
{\rm LT}(g_{\underline{m}})
\ \ \ \ \
\text{\rm for some}\ \
1\leqslant i_1<\cdots<i_{\underline{m}}\leqslant n.
\]
Again because ${\tt G}$ is reduced, each $g_l$ does not
contain any $x_{ i_1}, \dots, x_{ i_{ \underline{ m}}}$, aside
from its leading term $x_{ i_l}$. After relabelling
the $x_i$ if necessary, we can (and we shall) assume that 
$i_1 = n - \underline{ m} + 1$, \dots, $i_{\underline{ m}} =
n$. Then the $g_l$ write under a graphed form: 
\[
\aligned
g_l\big(x_1,\dots,x_{n-\underline{m}},
&\,\,
x_{n-\underline{m}+1},\dots,x_n\big)
=
x_l
-
g_l'(x_1,\dots,x_{n-\underline{m}})
\\
&
\ \ \ \ \ \ \ \ \ \ \ \ \ 
{\scriptstyle{(l\,=\,n\,-\,\underline{m}\,+\,1,\,\dots,\,n)}},
\endaligned
\]
for some linear forms $g_l'$ in only the $n - \underline{ m}$ first
variables $x_1, \dots, x_{n - \underline{m}}$.  
But then, since the vector subspace $F \subset E$ is as well
represented by the corresponding $\underline{ m}$ 
Cartesian linear equations $0 = x_l - g_l'(x_1, \dots,
x_{ \underline{ m}})$, for $l = 
n- \underline{ m} + 1, \dots, n$, it goes without saying that, 
in the notation of the proposition:
\[
m
:=
\dim_\K F 
=
n-\underline{m},
\]
so that we can replace $\underline{ m}$ by $n - m$ everywhere.
Furthermore, if we expand:
\[
g_l'(x_1,\dots,x_m)
=
\sum_{j=1}^m\,g_{lj}'\,x_j
\ \ \ \ \ \ \ \ \ \ \ \ \ 
{\scriptstyle{(l\,=\,m\,+\,1\,\cdots\,n)}}
\]
with some scalars $g_{lj}' \in \K$, it is clear that a certain basis
for $F$ which is naturally associated to the Cartesian linear
equations in question just consists of the
$m$ vectors obtained by setting one $x_j$ equal to $1$
and the others equal to $0$, for any choice of
$j = 1, \dots, m$, which yields the $m$ vectors:
\begin{equation}
\label{m-vectors}
{\sf e}_j
+
\sum_{l=m+1}^n\,g_{lj}'\,{\sf e}_l
\ \ \ \ \ \ \ \ \ \ \ \ \ {\scriptstyle{(j\,=\,1\,\cdots\,m)}}.
\end{equation}

On the other hand, the reduction of the auxiliary
bilinear form 
$\sum_{ i=1}^n\, x_i\, y_i$ to normal
form with respect to ${\tt G}$ then just means replacing $x_l$ by
$g_l' (x_1, \dots, x_m)$, for $l = m+1, \dots, n$, so that:
\[
\aligned
h_y(x)
=
{\rm NF}_{\tt G}
\Big(
{\textstyle{\sum_{i=1}^n}}\,x_iy_i
\Big)
&
=
\sum_{j=1}^n\,x_j\,y_j
+
\sum_{l=m+1}^n\,g_l'(x_1,\dots,x_m)\,y_l
\\
&
=
\sum_{j=1}^m\,x_j\,y_j
+
\sum_{l=m+1}^n\,\sum_{j=1}^m\,g_{lj}'\,x_j\,y_l
\\
&
=
\sum_{j=1}^m\,x_j
\bigg(
\underbrace{
y_j+\sum_{l=m+1}^n\,g_{lj}'\,y_l
}_{=:h_j(y)}
\bigg),
\endaligned
\] 
and from this last expression, one realizes that the 
$m$ vectors:
\[
{\sf h}_j
=
{\sf e}_j
+
\sum_{l=m+1}^n\,g_{lj}'\,{\sf e}_l
\ \ \ \ \ \ \ \ \ \ \ \ \ {\scriptstyle{(j\,=\,1\,\cdots\,m)}}
\]
associated to the obtained coefficients $h_j (y)$ of $h_y(x)$ with
respect to $x_1, \dots, x_m$ do indeed coincide with the $m = \dim F$
vectors~\thetag{ \ref{m-vectors}} which were seen to constitute a
basis for $F$ a moment ago.  The simultaneous proof of properties {\bf
(ii)}, {\bf (iii)}, {\bf (iv)} is therefore complete.
\endproof

\proof[Proof of Theorem~\ref{quotient-V-W}] After a permutation
of both the $\overline{ g}_j$ and
the variables $y_i$, we
can assume that the lexicographic ordering is just $y_n \prec \cdots
\prec y_2 \prec y_1$ and that the $q$ leading terms of 
the generators $\overline{
g}_1(y), \dots, \overline{ g}_q (y)$ of the Gr\"obner basis ${\tt B}_W$
are just $y_1, \dots, y_q$.
Since ${\tt B}_W$ is reduced, its $q$ elements
necessarily write under a graphed, linear form:
\[
{\tt B}_W
=
\Big\{
\underbrace{y_j
-
{\textstyle{\sum_{i=q+1}^{i=n}}}\,b_{j,i}\,y_i}_{
\overline{g}_j(y)}
\Big\}_{1\leqslant j\leqslant q},
\] 
for some scalars $b_{ {}_\bullet,
{}_\bullet} \in \K$.  Similarly, the $p$ elements 
$\overline{ f}_1 (y), \dots, \overline{ f}_p(y)$
of the
Gr\"obner basis ${\tt B}_V$ must also be of a certain graphed, linear
form.  Let $q' \leqslant q$ be the number of leading terms of elements
of ${\tt B}_V$ that are equal to one leading term $y_j$ with $1
\leqslant j \leqslant q$ appearing in the members of 
${\tt B}_W$.  Possibly after an
independent renumbering of both $y_1, \dots, y_q$ and $y_{ q+1},
\dots, y_n$, it follows that there is a decomposition of the
$y_i$-variables into four groups of variables:
\[
\big(
\underline{y_1,\dots,y_{q'}},
y_{q'+1},\dots,y_q,
\underline{y_{q+1},\dots,y_{p+q-q'}},
y_{p+q-q'+1},\dots,y_n
\big)
\]
such that the $p = q' + (p-q')$ elements of ${\tt B}_V$
do precisely have those leading monomials that are underlined and do
write under the following graphed form:
\[
\aligned
{\tt B}_V
&
=
\Big\{
y_{j'}
-
{\textstyle{\sum_{i=q'+1}^{i=q}}}\,a_{j',i}\,y_i
-
{\textstyle{\sum_{i=p+q-q'+1}^{i=n}}}\,a_{j',i}\,y_i
\Big\}_{1\leqslant j'\leqslant q'}
\bigcup
\\
&
\ \ \ \ \
\bigcup
\Big\{
y_l
-
{\textstyle{\sum_{i=q'+1}^{i=q}}}\,a_{l,i}\,y_i
-
{\textstyle{\sum_{i=p+q-q'+1}^{i=n}}}\,a_{l,i}\,y_i
\Big\}_{q+1\leqslant l\leqslant p+q-q'},
\endaligned
\]
for some scalars $a_{ {}_\bullet, {}_\bullet} \in \K$. 
However, all $a_{l,i}$ in the first sum of the second
line must necessarily be equal to $0$, because by
assumption, we have:
\[
y_l
\prec
y_{q'+1},\dots,y_q
\ \ \ \ \
\text{\rm for all}\ \
q+1\leqslant l\leqslant p+q-q',
\]
whence if some
$a_{ l,i}$ would
be nonzero, the number $q'$ defined above
would be larger. Thus, after simply erasing these $a_{ l,i}$, 
it remains:
\[
\aligned
{\tt B}_V
&
=
\Big\{
y_{j'}
-
{\textstyle{\sum_{i=q'+1}^{i=q}}}\,a_{j',i}\,y_i
-
{\textstyle{\sum_{i=p+q-q'+1}^{i=n}}}\,a_{j',i}\,y_i
\Big\}_{1\leqslant j'\leqslant q'}
\bigcup
\\
&
\ \ \ \ \
\bigcup
\Big\{
y_l
-
{\textstyle{\sum_{i=p+q-q'+1}^{i=n}}}\,a_{l,i}\,y_i
\Big\}_{q+1\leqslant l\leqslant p+q-q'}.
\endaligned
\]
But now, we remind the assumption $W \subset V$ which reads in terms
of ideals naturally as the constraint $\left< \, {\tt B}_W\, \right>
\subset \left<\, {\tt B}_V\, \right>$. Since all existing polynomials
are (degree-one) linear forms, each element $y_j - \sum_{ i=q+1}^{ i =
n}\, b_{j,i}\, y_i$ of ${\tt B}_W$ for $j = q'+1, \dots, q$ must in
particular be a certain linear combination of elements of ${\tt B}_V$
with scalar (degree-zero) coefficients. But all elements of ${\tt B}_V$
above are under a graphed form, with no such $y_j$ with $j = q'+1,
\dots, q$ appearing in either the $y_{ j'}$ or in the $y_l$ of ${\tt B}_V$,
from what we deduce $q' = q$, whence immediately:
\[
{\tt B}_V
=
\Big\{
y_j
-
{\textstyle{\sum_{i=p+1}^{i=n}}}\,
a_{j,i}\,y_i
\Big\}_{1\leqslant j\leqslant q}
\,\bigcup\,
\Big\{
y_l
-
{\textstyle{\sum_{i=p+1}^{i=n}}}\,
a_{l,i}\,y_i
\Big\}_{q+1\leqslant l\leqslant p}.
\]
Now that $q' = q$, the constraint $\left< \, {\tt B}_W\, \right>
\subset \left<\, {\tt B}_V\, \right>$ means that, for every
$j = 1, \dots, q$, there exist scalars $\lambda_{
j, j_1} \in \K$ and $\mu_{ j, l_1} \in \K$ such that one has:
\[
\aligned
y_j
-
{\textstyle{\sum_{i=q+1}^{i=n}}}\,
b_{j,i}\,y_i
&
\equiv
\sum_{j_1=1}^{j_1=q}\,
\lambda_{j,j_1}
\Big(
y_{j_1}
-
{\textstyle{\sum_{i=p+1}^{i=n}}}\,
a_{j_1,i}\,y_i
\Big)
+
\\
&
\ \ \ \ \
+
\sum_{l_1=q+1}^{l_1=p}\,
\mu_{j,l_1}
\Big(
y_{l_1}
-
{\textstyle{\sum_{i=p+1}^{i=n}}}\,
a_{l_1,i}\,y_i
\Big),
\endaligned
\]
identically in $\K [y_1, \dots, y_n]$. It
necessarily follows that $\lambda_{ j, j} = 1$ while
$\lambda_{ j, j_1} = 0$ for $j_1 \neq j$ and that:
\[
-\,b_{j,i}
=
\mu_{j,i}
\ \ \ \ \
\text{for}\ \
i=q+1,\dots,p.
\]
After simplifying terms which cancel out, there remain the 
$q$ equations:
\[
\aligned
-\,
{\textstyle{\sum_{i=p+1}^{i=n}}}\,
b_{j,i}\,y_i
&
\equiv
-\,
{\textstyle{\sum_{i=p+1}^{i=n}}}\,
a_{j,i}\,y_i
+
{\textstyle{\sum_{l_1=q+1}^{l_1=p}}}\,
{\textstyle{\sum_{i=p+1}^{i=n}}}\,
b_{j,l_1}\,a_{l_1,i}\,y_i
\\
&
\ \ \ \ \ \ \ \ \ \ \ \ \ \ \ 
{\scriptstyle{(j\,=\,1\,\cdots\,q)}},
\endaligned
\]
holding 
identically in $\K [y_1, \dots, y_n]$, and this yields by identification
of the coefficients of the $y_i$ in both sides:
\begin{equation}
\label{a-b-sum}
\aligned
b_{j,i}
&
=
a_{j,i}
-
{\textstyle{\sum_{l_1=q+1}^{l_1=p}}}\,
b_{j,l_1}\,a_{l_1,i}
\\
&
\ \ \ \ \ 
{\scriptstyle{(j\,=\,1\,\cdots\,q\,;\,\,
i\,=\,p\,+\,1\,\dots\,n)}}.
\endaligned
\end{equation}

On the other hand, reminding that computing the normal form with respect
to ${\tt B}_W$ just means replacing each $y_j$ by 
$\sum_{ i= q+1}^{ i=n}\, b_{ j,i}\, y_i$ for $j = 1, \dots, q$, 
we have:
\[
\aligned
\Big<\,
{\rm NF}_{{\tt B}_W}\big(\overline{f}\big)\colon
\overline{f}\in{\tt B}_V
\Big>
&
=
\Big<
\Big\{
{\textstyle{\sum_{i=q+1}^{i=n}}}\,b_{j,i}\,y_i
-
{\textstyle{\sum_{i=p+1}^{i=n}}}\,a_{j,i}\,y_i
\Big\}_{1\leqslant j\leqslant q},
\\
&
\ \ \ \ \ \ \ \ \
\Big\{
y_l
-
{\textstyle{\sum_{i=p+1}^{i=n}}}\,
a_{l,i}\,y_i
\Big\}_{q+1\leqslant l\leqslant p}
\Big>. 
\endaligned
\] 
In the first family, we use the relation~\thetag{ \ref{a-b-sum}}
obtained right above to replace the $b_{ j,i}$ for $1 \leqslant j
\leqslant q$ and for $p+1 \leqslant i \leqslant n$, which yields after
a cancellation:
\[
\aligned
\Big<\,
{\rm NF}_{{\tt B}_W}\big(\overline{f}\big)\colon
\overline{f}\in{\tt B}_V
\Big>
&
=
\Big<
\Big\{
{\textstyle{\sum_{i=q+1}^{i=p}}}\,b_{j,i}\,y_i
-
{\textstyle{\sum_{i=p+1}^{i=n}}}\,
{\textstyle{\sum_{l_1=q+1}^{l_1=p}}}\,
b_{j,l_1}\,a_{l_1,i}\,y_i
\Big\}_{1\leqslant j\leqslant q},
\\
&
\ \ \ \ \ \ \ \ \
\Big\{
y_l
-
{\textstyle{\sum_{i=p+1}^{i=n}}}\,
a_{l,i}\,y_i
\Big\}_{q+1\leqslant l\leqslant p}
\Big>. 
\endaligned
\]
But now, we observe that each element in the first family
belongs in fact already to the ideal generated by the members of the
second family, because the linear combination:
\[
{\textstyle{\sum_{l=q+1}^{l=p}}}\,
b_{j,l}\,
\Big(
y_l
-
{\textstyle{\sum_{i=p+1}^{i=n}}}\,
a_{l,i}\,y_i
\Big)
\]
identifies, after change of indices, 
to the $j$-th element of the first family. 
In conclusion, the ideal: 
\[
\aligned
\Big<\,
{\rm NF}_{{\tt B}_W}\big(\overline{f}\big)\colon
\overline{f}\in{\tt B}_V
\Big>
=
\Big<\,
\Big\{
y_l
-
{\textstyle{\sum_{i=p+1}^{i=n}}}\,
a_{l,i}\,y_i
\Big\}_{q+1\leqslant l\leqslant p}
\Big>
\endaligned
\]
is generated by exactly $p - q = \dim_\K V - \dim_\K W$ elements, with
are {\em de facto} in reduced Gr\"obner basis form 
for the lexicographic
ordering $\prec$, and the associated vectors:
\[
{\sf h}_l
=
{\sf e}_l
-
{\textstyle{\sum_{i=p+1}^{i=n}}}\,
a_{l,i}\,{\tt e}_i
\ \ \ \ \ \ \ \ \ \ \ \ \ 
{\scriptstyle{(l\,=\,q\,+\,1\,\cdots\,p)}}
\]
belong to $V$ by assumption (since vectors associated to 
elements of ${\tt B}_V$ belong to $V$) 
and are mutually linearly independent
modulo $W$, as one can easily realize thanks to the fact that $W$ is
graphed over $\K\, {\sf e}_1 \oplus \cdots \oplus \K\, {\sf e}_q$.
The proof of Theorem~\ref{quotient-V-W} is complete.
\endproof

\section{Description of the Algorithm based on Gr\"obner bases}
\label{Description} 

In this section we propose our new algorithm to compute the
cohomology spaces of Lie (super) algebras, based on
Proposition~\ref{extract-basis-from-Cartesian} and
Theorem~\ref{quotient-V-W}. This section includes also an example
which illustrates the behavior of this algorithm.

\begin{algorithm}[H]
\caption{{\sc ``LSAC''}}
\begin{algorithmic}

\REQUIRE 
$\left\{
\begin{array}{lll}
{\frak g}=\frak g_{\overline{0}}\oplus\frak g_{\overline{1}}
&:& 
\text{an $m$-dimensional Lie (super) algebra}
\\
{V}&:& 
\text{an $n$-dimensional $\frak g$-module}
\\
k&:& \text{the cohomology order}
\end{array}
\right.$

\smallskip\ENSURE{$H^k(\frak g,V)$};

\smallskip\STATE$\bullet$
${\rm Vars}
:=
\big\{\phi_{(i|j)_{r,k}}^l , 
\psi_{(i|j)_{r,k-1}}^l\big\}$;

\smallskip\STATE$\bullet$ 
$\prec\, :=$ a lexicographical ordering on 
$\mathbb{K}[{\rm Vars}]$ with
$\phi_{(i|j)_{r,k}}^l \prec \psi_{(i'|j')_{r,k-1}}^{l'}$;

\smallskip\STATE$\bullet$ 
${\sf Syst}_\phi(\mathcal{Z}^k) :=$ the set of equations
$(\partial^k\Phi)\big( ({\sf e}_i, {\sf o}_j)_{
s, k+1}\big) = 0$;

\smallskip\STATE$\bullet$ 
${\tt G}_{{\mathcal Z}^k} :=$ 
the reduced Gr\"obner basis of $\left< 
{\sf Syst}_\phi(\mathcal{Z}^k) \right>$ 
with respect to $\prec$;

\smallskip\STATE$\bullet$ 
${\sf Syst}_{\psi,\phi}(\mathcal{B}^k) :=$ 
the set of equations 
$\partial^{ k-1} \Psi \big( ({\sf e}_i, {\sf o}_j)_{ r, k} 
\big) = \Phi \big( ({\sf e}_i, {\sf o}_j)_{ r, k} \big)$;

\smallskip\STATE$\bullet$ 
${\tt G}_{{\mathcal B}^k} :=$ the reduced Gr\"obner basis of
$\left< {\sf Syst}_{\psi,\phi}(\mathcal{B}^k)
\right> \cap \mathbb{K}\big[\phi^{(i|j)_{r,k}}_l\big]$;

\smallskip\STATE$\bullet$ 
$\big\{\upsilon^{(i|j)_{r,k}}_l \big\}
:=$ the
auxiliary variables with the same order $\prec$ as for
$\big\{ \phi^{(i|j)_{r,k}}_l \big\}$;

\smallskip\STATE$\bullet$ 
${\sf BilinearForm} := 
\sum\,\phi^{ (i|j)_{ r,k}}_l\, \upsilon^{(i|j)_{r,k}}_l$
the auxiliary bilinear form in the two collections of variables 
$\phi^{\cdot}_{\cdot}$ and $\upsilon_\cdot^\cdot$; 

\smallskip\STATE$\bullet$ 
${\sf C}_{{\mathcal Z}^k}
:=
\Big\{
{\rm Coeff\!}_{\upsilon^{(i|j)_{r,k}}_l}
\Big({\rm NF}_{{\tt G}_{{\mathcal Z}^k}}
\big(
{\sf BilinearForm}
\big)
\Big)
\Big\}$;

\smallskip\STATE$\bullet$ 
${\tt Basis}({\mathcal Z}^k) :=$ the reduced Gr\"obner basis of
${\sf C}_{{\mathcal Z}^k}$ with
respect to $\prec$;

\smallskip\STATE$\bullet$ 
${\sf C}_{{\mathcal B}^k}
:=
\Big\{
{\rm Coeff\!}_{\upsilon^{(i|j)_{r,k}}_l}
\Big({\rm NF}_{{\tt G}_{{\mathcal B}^k}}
\big(
{\sf BilinearForm}
\big)
\Big)
\Big\}$;

\smallskip\STATE$\bullet$ 
${\tt Basis}({\mathcal B}^k) :=$ 
the reduced Gr\"obner basis of ${\sf C}_{{\mathcal B}^k}$ with
respect to $\prec$;

\medskip
\hspace{-0.5cm}
{\bf Return:}
${\tt Basis} \big( \mathcal{ Z}^k / 
\mathcal{ B}^k\big) :=$ the reduced Gr\"obner basis of:
\[
\left<\, 
{\rm NF}_{{\tt Basis}({\mathcal B}^k)}
\big(\vartheta\big) 
\colon 
\vartheta
\in 
{\tt Basis}({\mathcal Z}^k) \,\right>.
\]
\end{algorithmic}
\end{algorithm}

\begin{Example}
\label{MS-Example}
Let $\mathfrak{ h}$ be the $7$-dimensional
standard Lie algebra over $\mathbb{
Q}$ whose basis elements $\{{\sf l}_1,{\sf l}_2,{\sf d},{\sf t}_1,{\sf
t}_2,{\sf t}_3,{\sf r}\}$ enjoy the following commutator table
(\cite{ MS}):

\begin{center}
\begin{tabular} [t] { c | c c c c c c c } 
& ${\sf t}_1$ & ${\sf t}_2$ & ${\sf t}_3$ 
& ${\sf l}_1$ & ${\sf l}_2$ & ${\sf
r}$ & ${\sf d}$
\\
\hline ${\sf t}_1$ & $0$ & $0$ & $0$ 
& $-\,{\sf t}_2$ & $-{\sf t}_3$ & $0$ & $2{\sf t}_1$
\\
${\sf t}_2$ & $*$ & $0$ & $0$ & $0$ & $0$ & ${\sf t}_3$ & $-3\,{\sf t}_2$
\\
${\sf t}_3$ & $*$ & $*$ & $0$ & $0$ & $0$ 
& $-\,{\sf t}_2$ & $-3\,{\sf t}_3$
\\
${\sf l}_1$ & $*$ & $*$ & $*$ & $0$ & ${\sf t}_1$ 
& ${\sf l}_2$ & $-{\sf l}_1$
\\
${\sf l}_2$ & $*$ & $*$ & $*$ & $*$ & $0$ & $-{\sf l}_1$ & $-{\sf l}_2$
\\
${\sf r}$ & $*$ & $*$ & $*$ & $*$ & $*$ & $0$ & $0$
\\
${\sf d}$ & $*$ & $*$ & $*$ & $*$ & $*$ & $*$ & $0$
\end{tabular}
\end{center}

\noindent
and let $\frak g$ be the Lie subalgebra of $\mathfrak{ h}$ which is
generated by $\{{\sf l}_1,{\sf
l}_2,{\sf t}_1,{\sf t}_2,{\sf t}_3\}$. We would like to compute the
fourth cohomology space $H^4 ( \mathfrak{ g}, \mathfrak{
h})$. Applying the algorithm, a computer yields the reduced Gr\"obner
basis: 
\[
\footnotesize\aligned
{\tt G}_{\mathcal{Z}^4}
&
=
\Big\{
\phi_{l_1, t_1, t_2, t_3}^r
-
\phi_{l_2, t_1, t_2, t_3}^d,\,\,\,
\phi_{l_1, t_1, t_2, t_3}^d
+
\phi_{l_2, t_1, t_2, t_3}^r,\,\,\,
2\phi_{l_1, l_2, t_2, t_3}^d
-
\phi_{l_1, t_1, t_2, t_3}^{l_1}
-
\phi_{l_2,t_1,t_2,t_3}^{l_2},
\\
& 
\ \ \ \ \ \ \ \ 
\phi_{l_1,l_2,t_1,t_2}^r-3\phi_{l_1,
l_2, t_1, t_3}^d
+
\phi_{l_1, l_2, t_2, t_3}^{l_1}-
\phi_{l_2, t_1, t_2, t_3}^{t_1},
\\
&
\ \ \ \ \ \ \ \ 
3\phi_{l_1,l_2, t_1, t_2}^d
+
\phi_{l_1, l_2, t_1, t_3}^r
+
\phi_{l_1, l_2, t_2, t_3}^{l_2}
+
\phi_{l_1, t_1, t_2,t_3}^{t_1}
\Big\},
\endaligned
\]
together with:
\[
\footnotesize\aligned
{\tt G}_{\mathcal{B}^4}
&
=
\Big\{
\phi_{l_2, t_1, t_2, t_3}^r,\,\,\,
\phi_{l_2, t_1, t_2, t_3}^d,\,\,\,
\phi_{l_1, t_1, t_2, t_3}^r,\,\,\, 
\phi_{l_1, t_1, t_2, t_3}^d,\,\,\,
\phi_{l_1, t_1, t_2, t_3}^{ l_2}
+
\phi_{l_2, t_1, t_2, t_3}^{l_1},\,\,\,
\\
& 
\ \ \ \ \ \ \ \ 
-\,\phi_{l_2,t_1,t_2,t_3}^{l_2}
+
\phi_{l_1,t_1,t_2,t_3}^{l_1},\,\,\,
-\,\phi_{l_2, t_1, t_2, t_3}^{l_1}
+
\phi_{l_1, l_2, t_2, t_3}^r,\,\,\, 
-
\phi_{l_2, t_1, t_2,t_3}^{ l_2}
+
\phi_{l_1, l_2, t_2, t_3}^{ d},\,\,\, 
\\
&
\ \ \ \ \ \ \ \ 
\phi_{l_1, l_2, t_1,t_2}^r
-
3\phi_{l_1, l_2, t_1, t_3}^d
+
\phi_{l_1, l_2, t_2, t_3}^{l_1}
-
\phi_{l_2, t_1, t_2, t_3}^{ t_1},
\\
&
\ \ \ \ \ \ \ \ 
3\phi_{l_1, l_2, t_1,t_2}^d
+
\phi_{l_1, l_2, t_1, t_3}^r
+
\phi_{l_1, l_2, t_2, t_3}^{l_2}
+
\phi_{l_1, t_1, t_2, t_3}^{ t_1}
\Big\}.
\endaligned
\]
Next, relabelling  the variables $\phi_\cdot^\cdot$ and
$\upsilon_\cdot^\cdot$ by $x_1, \dots, x_{ 35}$ and
$y_1, \dots, y_{ 35}$, we obtain:
\[
\footnotesize
\aligned 
{\tt Basis}({\mathcal Z}^4)
& 
=
\big\{
x_{34},\,\, 
x_{33},\,\, x_{29},\,\, x_{28}+x_{31},\,\,
x_{27},\,\, x_{26},\,\, x_{24}-x_{35},\,\, x_{23},\,\, 
-x_{30}+x_{22},\,\, x_{21},\,\, x_{20},\,\, x_{19},\,\,
\\
&
\ \ \ \ \ \ \ \ 
x_{18},\,\, x_{17}+2x_{30},\,\, -x_{25}+x_{16},\,\, x_{15}+x_{32},\,\, -x_{25}+x_{14},\,\, x_{13},\,\, x_{12},\,\, x_{11},\,\,
-3x_{32}+x_{10},\,\,
\\
&
\ \ \ \ \ \ \ \ 
x_{9},\,\, x_{8},\,\, x_{32}+x_{7},\,\, x_{6},\,\, x_{5},\,\, x_{4},\,\, -3x_{25}+x_{3},\,\, x_{2},\,\, x_{1}
\big\},
\endaligned
\]
\[\footnotesize
\aligned
{\tt Basis}({\mathcal B}^4) 
&
=
\big\{
x_{34},\,\, x_{33},\,\, x_{27},\,\, x_{26},\,\,
-x_{23}+x_{21}+x_{29},\,\, x_{20},\,\,
x_{19},\,\, x_{18},\,\, x_{22}+x_{17}+x_{30},\,\, -x_{25}+x_{16},\,\,
\\
&
\ \ \ \ \ \ \ \ 
x_{15}+x_{32},\,\, -x_{25}+x_{14},\,\, x_{13},\,\, x_{12},\,\, x_{11},\,\,
-3x_{32}+x_{10},\,\, x_{9},\,\, x_{8},\,\, x_{32}+x_{7},\,\, x_{6},\,\, x_{5},\,\, x_{4},\,\, 
\\
&
\ \ \ \ \ \ \ \ 
-3x_{25}+x_{3},\,\, x_{2},\,\, x_{1}
\big\},
\endaligned
\]
of cardinalities $30$ and $25$, respectively. The last
step provides a basis of $5 = 30 - 25$ vectors for 
$\mathcal{ Z}^4 / \mathcal{ B}^4$ represented by 
means of the following $5$ associated linear forms:
\[
{\tt Basis}
\big(\mathcal{Z}^4\big/{\mathcal B}^4\big) 
=
\big\{
x_{29},\,\,
x_{28}+x_{31},\,\, 
x_{24}-x_{35},\,\, 
x_{23},\,\, 
-x_{30}+x_{22}
\big\},
\]
and coming back to the original notation, this
corresponds to:
\[
\aligned
{\tt Basis}
\big(\mathcal{Z}^4\big/{\mathcal B}^4\big) 
= 
\Big\{ 
&
{\sf l}_2^*\wedge{\sf t}_1^*\wedge{\sf t}_2^*\wedge{\sf t}_3^*
\otimes{\sf l}_1,\,\,\,
{\sf l}_1^*\wedge{\sf t}_1^*\wedge{\sf t}_2^*\wedge{\sf t}_3^*
\otimes{\sf r}
+
{\sf l}_2^*\wedge{\sf t}_1^*\wedge{\sf t}_2^*\wedge{\sf t}_3^*
\otimes{\sf d},
\\
& 
{\sf l}_1^*\wedge{\sf t}_1^*\wedge{\sf t}_2^*\wedge{\sf t}_3^*
\otimes{\sf l}_2,\,\,\,
{\sf l}_1^*\wedge{\sf t}_1^*\wedge{\sf t}_2^*\wedge{\sf t}_3^*
\otimes{\sf d}
-
{\sf l}_2^*\wedge{\sf t}_1^*\wedge{\sf t}_2^*\wedge{\sf t}_3^*
\otimes{\sf r},
\\
&
{\sf l}_1^*\wedge{\sf t}_1^*\wedge{\sf t}_2^*\wedge{\sf t}_3^*
\otimes{\sf l}_1
-
{\sf l}_2^*\wedge{\sf t}_1^*\wedge{\sf t}_2^*\wedge{\sf t}_3^*
\otimes{\sf l}_2
\Big\}.
\endaligned
\]
\end{Example}

\section{Improvement of the Algorithm when Cohomology Spaces Split}
\label{Improved}

As we saw, the two collections of Cartesian linear equations ${\sf
Syst}_\phi (\mathcal{ Z}^k)$ and ${\sf Syst}_\phi (\mathcal{ Z}^k)$
have an essential r\^ole in the process, and if the number of
variables in them increases, one can expect that the complexity of
computations will increases too.  Here, in the case of standard Lie
algebras $\mathfrak{ g} \subset \mathfrak{ h} = V$, one further aim
could to set up a refined algorithm which inspects whether these
equations split up into a collection of sub-equations each of which
involves a smaller number of variables. However, this kind
of problem lies a bit outside the scope of the present
article, closer to plain searching-and-listing algorithmic procedures, 
because it amounts to read, by means of
a computer, some two given systems of linear equations in some variables
$(x_1, \dots, x_n)$ and to pick up step by step the appearing
nonzero $\lambda_i\, x_i$ until one gathers pairs
of collections of equations which involve only
a {\em subset} of variables, all subsets being pairwise distinct.

Nevertheless, the circumstance of spitting up
naturally occurs for instance when the Lie algebras $\mathfrak{ g}$
and $\mathfrak{ h}$ are {\sl graded}
at the beginning, in the sense of Tanaka (\cite{
Tanaka, AMS}), namely when one has
decompositions into direct sums of $\K$-vector subspaces:
\[
\aligned
\mathfrak{h}
&
=
\mathfrak{h}_{-a}
\oplus\cdots\oplus
\mathfrak{h}_{-1}
\oplus
\mathfrak{h}_0
\oplus
\mathfrak{h}_1
\oplus\cdots\oplus
\mathfrak{h}_b
\\
\mathfrak{g}
&
=
\mathfrak{h}_{-a}
\oplus\cdots\oplus
\mathfrak{h}_{-1},
\endaligned
\] 
where $a \geqslant 1$ and $b \geqslant 0$ are certain
two integers, with
the property that:
\[
\big[\mathfrak{h}_{\ell_1},\,\mathfrak{h}_{\ell_2}\big]
\subset 
\mathfrak{h}_{\ell_1+\ell_2},
\] 
for all $\ell_1, \ell_2 \in \mathbb{ Z}$, after prolonging trivially
$\mathfrak{ h}_\ell 
:= \{ 0\}$ for either $\ell \leqslant - a - 1$ or $\ell
\geqslant b +1$. Then each space of $k$-cochains
$\mathcal{ C}^k (\mathfrak{
g}, \mathfrak{ h})$ naturally splits up as a direct sum of so-called
{\sl homogeneous $k$-cochains} as follows: a $k$-cochain
$\Phi \in \mathcal{C}^k ( \mathfrak{ g}, \mathfrak{ h})$ is said to
be {\sl of homogeneity} a certain integer $h \in \mathbb{ Z}$ whenever
for any $k$ vectors:
\[
{\sf z}_{i_1}\in\mathfrak{h}_{\ell_1}, 
\ldots\ldots, 
{\sf z}_{i_k}\in\mathfrak{h}_{\ell_k}
\]
belonging to certain arbitrary but determined $\mathfrak{ h}$-components, 
its value: 
\[
\Phi({\sf z}_{i_1},\dots,{\sf z}_{i_k})
\in
\mathfrak{h}_{\ell_1+\cdots+\ell_k+h}
\]
belongs to the $(\ell_1 + \cdots + \ell_k + h)$-th component of
$\mathfrak{ h}$. Then one easily convinces oneself 
({\em see} also~\cite{ Goze})
that any 
$k$-cochain $\Phi \in \mathcal{ C}^k
( \mathfrak{ g}, \mathfrak{ h})$ 
splits up as a direct sum of $k$-cochains of
fixed homogeneity:
\[
\Phi
=
\cdots+
\Phi^{[h-1]}
+
\Phi^{[h]}
+
\Phi^{[h+1]}
+\cdots,
\]
where we denote the completely $h$-homogeneous component of
$\Phi$ just by $\Phi^{ [h]}$. In other words:
\[
\mathcal{C}^k(\mathfrak{g},\mathfrak{h})
=
\bigoplus_{h\in\mathbb{Z}}\,
\mathcal{C}_{[h]}^k(\mathfrak{g},\mathfrak{h}),
\]
where of course the spaces $\mathcal{ C}_{ [h]}^k ( \mathfrak{ g},
\mathfrak{ h})$ reduce to $\{ 0\}$ for all large $\vert h
\vert$. Furthermore, applying the definition~\thetag{
\ref{standard-cochain}}, one verifies the important fact that
$\partial^k$ respects homogeneity for all $k = 0, 1, \dots, n$,
that is to say, for any $h \in \mathbb{Z}$, one has: 
\[
\partial^k\big(\mathcal{C}_{[h]}^k\big) 
\subset 
\mathcal{C}_{[h]}^{ k+1},
\]
whence the complex~\thetag{ \ref{complex-partial}} splits up as a
direct sum of complexes:
\[
0
\overset{\partial_{[h]}^0}{\longrightarrow}
\mathcal{C}^1
\overset{\partial_{[h]}^1}{\longrightarrow}
\mathcal{C}^2
\overset{\partial_{[h]}^2}{\longrightarrow}
\cdots
\overset{\partial_{[h]}^{m-2}}{\longrightarrow}
\mathcal{C}^{m-1}
\overset{\partial_{[h]}^{m-1}}{\longrightarrow}
\mathcal{C}^m
\overset{\partial_{[h]}^m}{\longrightarrow}
0
\]
indexed by $h \in \mathbb{Z}$, where $\partial_{ [h]}^k$ naturally denotes the
restriction:
\[
\partial_{[h]}^k
:=
\partial^k\big\vert_{\mathcal{C}_{[h]}^k}
\colon
\mathcal{C}_{[h]}^k
\longrightarrow
\mathcal{C}_{[h]}^{k+1}.
\]
Consequently, one may introduce the spaces of
{\sl $h$-homogeneous cocycles} of order $k$: 
\[
\mathcal{Z}_{[h]}^k(\mathfrak{g},\mathfrak{h}\big)
:=
{\rm ker}\big(\partial_{[h]}^k
\colon
\mathcal{C}_{[h]}^k\to\mathcal{C}_{[h]}^{k+1}\big),
\]
together with the spaces of 
{\sl $h$-homogeneous coboundaries} of order $k$:
\[
\mathcal{B}_{[h]}^k(\mathfrak{g},\mathfrak{h}\big)
:=
{\rm im}\big(\partial_{[h]}^{k-1}
\colon
\mathcal{C}_{[h]}^{k-1}\to\mathcal{C}_{[h]}^k\big).
\]
The computation of the
$h$-homogeneous $k$-th cohomology spaces:
\[
H_{[h]}^k\big(\mathfrak{g},\mathfrak{h}\big)
:=
\frac{\mathcal{Z}_{[h]}^k(\mathfrak{g},\mathfrak{h}\big)}{
\mathcal{B}_{[h]}^k(\mathfrak{g},\mathfrak{h}\big)}
\]
then requires to deal with vector (sub)spaces of smaller dimensions
and enables one to reconstitute the complete
cohomology space: 
\[
H^k(\mathfrak{g},\mathfrak{g})
=
\bigoplus_{h\in\mathbb{Z}}\,
H_{[h]}^k(\mathfrak{g},\mathfrak{g}).
\]

\begin{Example}
\label{AMS-Example}
Let $\mathfrak{ h}$ be the
$8$-dimensional Lie algebra over $\mathbb{ Q}$ whose basis
elements $\{ {\sf t}, {\sf h}_1, {\sf h}_2, {\sf r}, {\sf d},
{\sf i}_1, {\sf i}_2, {\sf j} \}$
enjoy the following commutator table:

\medskip
\begin{center}
\begin{tabular} [t] { l | l l l l l l l l } 
& ${\sf t}$ & ${\sf h}_1$ & ${\sf h}_2$ & ${\sf d}$ & ${\sf r}$ & ${\sf i}_1$ &
${\sf i}_2$ & ${\sf j}$
\\
\hline ${\sf t}$ & $0$ & $0$ & $0$ & $2\,{\sf t}$ & $0$ & ${\sf h}_1$ & ${\sf h}_2$ & ${\sf d}$
\\
${\sf h}_1$ & $*$ & $0$ & $4\,{\sf t}$ & ${\sf h}_1$ & ${\sf h}_2$ & $6\,{\sf r}$ & $2\,{\sf d}$ &
${\sf i}_1$
\\
${\sf h}_2$ & $*$ & $*$ & $0$ & ${\sf h}_2$ & $-{\sf h}_1$ & $-2\,{\sf d}$ & $6\,{\sf r}$ & ${\sf
i}_2$
\\
${\sf d}$ & $*$ & $*$ & $*$ & $0$ & $0$ & ${\sf i}_1$ & ${\sf i}_2$ & $2\,{\sf j}$
\\
${\sf r}$ & $*$ & $*$ & $*$ & $*$ & $0$ & $-{\sf i}_2$ & ${\sf i}_1$ & $0$
\\
${\sf i}_1$ & $*$ & $*$ & $*$ & $*$ & $*$ & $0$ & $4\,{\sf j}$ & $0$
\\
${\sf i}_2$ & $*$ & $*$ & $*$ & $*$ & $*$ & $*$ & $0$ & $0$
\\
${\sf j}$ & $*$ & $*$ & $*$ & $*$ & $*$ & $*$ & $*$ & $0$
\end{tabular}
\end{center}

\noindent
and let $\mathfrak{ g}$ be the Lie subalgebra of $\mathfrak{ h}$ which
is generated by ${\sf t}, {\sf h}_1, {\sf h}_2$, {\em see}~\cite{ AMS}
for application to the differential study of Cartan connection in
local Cauchy-Riemann geometry. We want to compute $H^2 ( \mathfrak{
g}, \mathfrak{ h})$. The geometry provides a natural graduation:
\[
\mathfrak{h}
=
\underbrace{
\mathfrak{h}_{-2}
\oplus 
\mathfrak{h}_{-1}}_{\mathfrak{g}}
\oplus
\mathfrak{h}_0
\oplus 
\mathfrak{h}_1
\oplus
\mathfrak{h}_2
\]
where:
\[
\mathfrak{h}_{-2}
=
\mathbb{R}\,{\sf t},\ \ \
\mathfrak{h}_{-1}
=
\mathbb{R}\,{\sf h}_1\oplus\mathbb{R}\,{\sf h}_2,\ \ \
\mathfrak{h}_0
=
\mathbb{R}\,{\sf d}\oplus\mathbb{R}\,{\sf r},\ \ \
\mathfrak{h}_1
=
\mathbb{R}\,{\sf i}_1\oplus\mathbb{R}\,{\sf i}_2,\ \ \
\mathfrak{h}_2
=
\mathbb{R}\,{\sf j},
\]
and one verifies that the commutator table written above
respects this graduation. A general $2$-cochain 
$\Phi \in \Lambda^2 \mathfrak{ g}^* \otimes \mathfrak{ h}$
writes under the
form:
\[
\scriptsize
\aligned
\Phi
&
=
\phi_t^{h_1h_2}\,{\sf h}_1^*\!\!\wedge{\sf h}_2^*\otimes{\sf t}
+
\ \ \ \ \
\fbox{\tiny 0}
\\
\fbox{\tiny 1}
&
\ \ \ \ \
+
\phi_t^{th_1}\,{\sf t}^*\!\!\wedge{\sf h}_1^*\otimes{\sf t}
+
\phi_t^{th_2}\,{\sf t}^*\!\!\wedge{\sf h}_2^*\otimes{\sf t}
+
\phi_{h_1}^{h_1h_2}\,{\sf h}_1^*\!\!\wedge{\sf h}_2^*\otimes{\sf t}
+
\phi_{h_2}^{h_1h_2}\,{\sf h}_1^*\!\!\wedge{\sf h}_2^*\otimes{\sf h}_2
+
\\
\fbox{\tiny 2}
&
\ \ \ \ \
+
\phi_{h_1}^{th_1}\,{\sf t}^*\!\!\wedge{\sf h}_1^*\otimes{\sf h}_1
+
\phi_{h_2}^{th_1}\,{\sf t}^*\!\!\wedge{\sf h}_1^*\otimes{\sf h}_2
+
\phi_{h_1}^{th_2}\,{\sf t}^*\!\!\wedge{\sf h}_2^*\otimes{\sf h}_1
+
\phi_{h_2}^{th_2}\,{\sf t}^*\!\!\wedge{\sf h}_2^*\otimes{\sf h}_2
+
\\
&
\ \ \ \ \ \ \ \ \ \ \ \ \ \ \ \ \ \ \ \ \ \ \ \ \ \ \ \ \ \ \ \ \ \
\ \ \ \ \ \ \ \ \ \ \ \ \ \ \ \ \ \ \ \ \ \ \ \ \ \ \ \ \ \ \ \ \ \ 
\ \ \ \
+
\phi_d^{h_1h_2}\,{\sf h}_1^*\!\!\wedge{\sf h}_2^*\otimes{\sf d}
+
\phi_r^{h_1h_2}\,{\sf h}_1^*\!\!\wedge{\sf h}_2^*\otimes{\sf r}
+
\\
\fbox{\tiny 3}
&
\ \ \ \ \
+
\phi_d^{th_1}\,{\sf t}^*\!\!\wedge{\sf h}_1^*\otimes{\sf d}
+
\phi_r^{th_1}\,{\sf t}^*\!\!\wedge{\sf h}_1^*\otimes{\sf r}
+
\phi_d^{th_2}\,{\sf t}^*\!\!\wedge{\sf h}_2^*\otimes{\sf d}
+
\phi_r^{th_2}\,{\sf t}^*\!\!\wedge{\sf h}_2^*\otimes{\sf r}
\\
&
\ \ \ \ \ \ \ \ \ \ \ \ \ \ \ \ \ \ \ \ \ \ \ \ \ \ \ \ \ \ \ \ \ \
\ \ \ \ \ \ \ \ \ \ \ \ \ \ \ \ \ \ \ \ \ \ \ \ \ \ \ \ \ \ \ \ 
+
\phi_{i_1}^{h_1h_2}\,{\sf h}_1^*\!\!\wedge{\sf h_2}^*\otimes{\sf i}_1
+
\phi_{i_2}^{h_1h_2}\,{\sf h}_1^*\!\!\wedge{\sf h}_2^*\otimes{\sf i}_2
+
\\
\fbox{\tiny 4}
&
\ \ \ \ \
+
\phi_{i_1}^{th_1}\,{\sf t}^*\!\!\wedge{\sf h}_1^*\otimes{\sf i}_1
+
\phi_{i_2}^{th_1}\,{\sf t}^*\!\!\wedge{\sf h}_1^*\otimes{\sf i}_2
+
\phi_{i_1}^{th_2}\,{\sf t}^*\!\!\wedge{\sf h}_2^*\otimes{\sf i}_1
+
\phi_{i_2}^{th_2}\,{\sf t}^*\!\!\wedge{\sf h}_2^*\otimes{\sf i}_2
\\
&
\ \ \ \ \ \ \ \ \ \ \ \ \ \ \ \ \ \ \ \ \ \ \ \ \ \ \ \ \ \ \ \ \ \
\ \ \ \ \ \ \ \ \ \ \ \ \ \ \ \ \ \ \ \ \ \ \ \ \ \ \ \ \ \ \ \ \ \ 
\ \ \ \ \ \ \ \ \ \ \ \ \ \ \ \ \ \ \ \ \ \ \ \ \ \ \ \ \ \ \ \ \, 
+
\phi_j^{h_1h_2}\,{\sf h}_1^*\!\!\wedge{\sf h}_2^*\otimes{\sf j}
+
\\
\fbox{\tiny 5}
&
\ \ \ \ \
+
\phi_j^{th_1}\,{\sf t}^*\!\!\wedge{\sf h}_1^*\otimes{\sf j}
+
\phi_j^{th_2}\,{\sf t}^*\!\!\wedge{\sf h}_1^*\otimes{\sf j},
\endaligned
\]
where framed numbers denote homogeneity of their lines. 
After computations, 
a $2$-cochain $\Phi$ is a $2$-cocycle if and only if its 24
coefficients satisfy the following seven linear equations, ordered
line by line by
increasing homogeneity:
\[
\footnotesize
\aligned
&
\fbox{\tiny 2}
\ \ \ \ \ \ \ \ \ \
0
=
2\phi_d^{h_1h_2}-4\phi_{h_2}^{th_2}-4\phi_{h_1}^{th_1},
\\
&
\fbox{\tiny 3}
\ \ \ \ \ \ \ \ \ \
0
=
\phi_{i_1}^{h_1h_2}-\phi_d^{th_2}-\phi_r^{th_1},
\ \ \ \ \ \ \ \ \ \ \ \ \ \
0
=
\phi_{i_2}^{h_1h_2}-\phi_r^{th_2}+\phi_d^{th_1},
\\
&
\fbox{\tiny 4}
\ \ \ \ \ \ \ \ \ \
0
=
\phi_j^{h_1h_2}-2\phi_{i_2}^{th_2}-2\phi_{i_1}^{th_1},
\ \ \ \ \ \ \ \ \ \
0
=
-6\phi_{i_1}^{th_2}+6\phi_{i_2}^{th_1},
\\
&
\fbox{\tiny 5}
\ \ \ \ \ \ \ \ \ \
0
=
-\phi_j^{th_2},
\ \ \ \ \ \ \ \ \ \
0
=
\phi_j^{th_1}.
\endaligned
\]
Next, a general $1$-cochain $\Psi \in \Lambda^1 \mathfrak{ g}^*
\otimes \mathfrak{ h}$ writes under the form: 
\[
\footnotesize
\aligned
\Psi
&
=
\psi_t^{h_1}\,{\sf h}_1^*\otimes{\sf t}
+
\psi_t^{h_2}\,{\sf h}_2^*\otimes{\sf t}
+
\ \ \ \ \ \ \ \ \ \ \
\fbox{\tiny -1}
\\
\fbox{\tiny 0}
&
\ \ \ \ \
+
\psi_t^t\,{\sf t}^*\otimes{\sf t}
+
\psi_{h_1}^{h_1}\,{\sf h}_1^*\otimes{\sf h}_1
+
\psi_{h_2}^{h_1}\,{\sf h}_1^*\otimes{\sf h}_2
+
\psi_{h_1}^{h_2}\,{\sf h}_2^*\otimes{\sf h}_1
+
\psi_{h_2}^{h_2}\,{\sf h}_2^*\otimes{\sf h}_2
+
\\
\fbox{\tiny 1}
&
\ \ \ \ \
+
\psi_{h_1}^t\,{\sf t}^*\otimes{\sf h}_1
+
\psi_{h_2}^t\,{\sf t}^*\otimes{\sf h}_2
+
\psi_d^{h_1}\,{\sf h}_1^*\otimes{\sf d}
+
\psi_r^{h_1}\,{\sf h}_1^*\otimes{\sf r}
+
\psi_d^{h_2}\,{\sf h}_2^*\otimes{\sf d}
+
\psi_r^{h_2}\,{\sf h}_2^*\,\otimes{\sf r}
+
\\
\fbox{\tiny 2}
&
\ \ \ \ \
+
\psi_d^t\,{\sf t}^*\otimes{\sf d}
+
\psi_r^t\,{\sf t}^*\otimes{\sf r}
+
\psi_{i_1}^{h_1}\,{\sf h}_1^*\otimes{\sf i}_1
+
\psi_{i_2}^{h_1}\,{\sf h}_1^*\otimes{\sf i}_2
+
\psi_{i_1}^{h_2}\,{\sf h}_2^*\otimes{\sf i}_1
+
\psi_{i_2}^{h_2}\,{\sf h}_2^*\otimes{\sf i}_2
+
\\
\fbox{\tiny 3}
&
\ \ \ \ \
+
\psi_{i_1}^t\,{\sf t}^*\otimes{\sf i}_1
+
\psi_{i_2}^t\,{\sf t}^*\otimes{\sf i}_2
+
\psi_j^{h_1}\,{\sf h}_1^*\otimes{\sf j}
+
\psi_j^{h_2}\,{\sf h}_2^*\otimes{\sf j}
+
\\
\fbox{\tiny 4}
&
\ \ \ \ \
+
\psi_j^t\,{\sf t}^*\otimes{\sf j}.
\endaligned
\]
The condition that $\Phi = \partial^1 \Psi$ then reads
in homogeneous-decomposed form:
\[
\footnotesize
\aligned
\fbox{\tiny 1}\ \ \ \ \
\phi_t^{th_1}
&
=
2\psi_d^{h_1}-4\psi_{h_2}^t
\\
\fbox{\tiny 2}\ \ \ \ \
\phi_{h_1}^{th_1}
&
=
\psi_{i_1}^{h_1}-\psi_d^t
\\
\fbox{\tiny 2}\ \ \ \ \
\phi_{h_2}^{th_1}
&
=
\psi_{i_2}^{h_1}-\psi_r^t
\\
\fbox{\tiny 3}\ \ \ \ \
\phi_d^{th_1}
&
=
\psi_j^{h_1}-2\psi_{i_2}^t
\\
\fbox{\tiny 3}\ \ \ \ \
\phi_r^{th_1}
&
=
-6\psi_{i_1}^t
\\
\fbox{\tiny 4}\ \ \ \ \
\phi_{i_1}^{th_1}
&
=
-\psi_j^t
\\
\fbox{\tiny 4}\ \ \ \ \
\phi_{i_2}^{th_1}
&
=
0
\\
\fbox{\tiny 5}\ \ \ \ \
\phi_j^{th_1}
&
=
0
\endaligned
\ \ \ \ \
\aligned
\fbox{\tiny 1}\ \ \ \ \
\phi_t^{th_2}
&
=
2\psi_d^{h_2}+4\psi_{h_1}^t
\\
\fbox{\tiny 2}\ \ \ \ \
\phi_{h_1}^{th_2}
&
=
\psi_{i_1}^{h_2}+\psi_r^t
\\
\fbox{\tiny 2}\ \ \ \ \
\phi_{h_2}^{th_2}
&
=
\psi_{i_2}^{h_2}-\psi_d^t
\\
\fbox{\tiny 3}\ \ \ \ \
\phi_d^{th_2}
&
=
\psi_j^{h_2}+2\psi_{i_1}^t
\\
\fbox{\tiny 3}\ \ \ \ \
\phi_r^{th_2}
&
=
-6\psi_{i_2}^t
\\
\fbox{\tiny 4}\ \ \ \ \
\phi_{i_1}^{th_2}
&
=
0
\\
\fbox{\tiny 4}\ \ \ \ \
\phi_{i_2}^{th_2}
&
=
-\psi_j^t
\\
\fbox{\tiny 5}\ \ \ \ \
\phi_j^{th_2}
&
=
0
\endaligned
\ \ \ \ \
\aligned
\fbox{\tiny 0}\ \ \ \ \
\phi_t^{h_1h_2}
&
=
4\psi_{h_2}^{h_2}+4\psi_{h_1}^{h_1}-4\psi_t^t
\\
\fbox{\tiny 1}\ \ \ \ \
\phi_{h_1}^{h_1h_2}
&
=
\psi_d^{h_2}+\psi_r^{h_1}-4\psi_{h_1}^t
\\
\fbox{\tiny 1}\ \ \ \ \
\phi_{h_2}^{h_1h_2}
&
=
\psi_r^{h_2}-\psi_d^{h_1}+4\psi_{h_2}^t
\\
\fbox{\tiny 2}\ \ \ \ \
\phi_d^{h_1h_2}
&
=
2\psi_{i_2}^{h_2}+2\psi_{i_1}^{h_1}-4\psi_d^t
\\
\fbox{\tiny 2}\ \ \ \ \
\phi_r^{h_1h_2}
&
=
6\psi_{i_1}^{h_2}-6\psi_{i_2}^{h_1}-4\psi_r^t
\\
\fbox{\tiny 3}\ \ \ \ \
\phi_{i_1}^{h_1h_2}
&
=
\psi_j^{h_2}-4\psi_{i_1}^t
\\
\fbox{\tiny 3}\ \ \ \ \
\phi_{i_2}^{h_1h_2}
&
=
-\psi_j^{h_1}-4\psi_{i_2}^t
\\
\fbox{\tiny 4}\ \ \ \ \
\phi_j^{h_1h_2}
&
=
-4\psi_j^t.
\endaligned
\]
One can then apply our algorithm to each subcollection of 
equations for every fixed homogeneity, and find that
$H^2 ( \mathfrak{ g}, \mathfrak{ h})$ is $2$-dimensional, 
generated by:
\[
\boxed{
\aligned
&
{\sf t}^*\wedge{\sf h}_2^*
\otimes
{\sf i}_2
-
2{\sf h}_1^*\wedge
{\sf h}_2^*
\otimes
{\sf j}
\\
\text{\rm and:}
\ \ \ \ \
&
{\sf t}^*\wedge{\sf h}_2^*
\otimes
{\sf i}_1
-
{\sf t}^*\wedge{\sf h}_1^*
\otimes
{\sf i}_2,
\endaligned}
\]
with the further observation that all cohomologies are zero
except in homogeneity $4$:
\begin{center}
\label{dimensional-cohomologies}
\begin{tabular}{|c|c|c|c|c|}
\hline\vspace{-10pt} &&&\\
\text{\rm Homogeneity} & $\dim\mathcal{C}^2$ & $\dim\mathcal{Z}^2$
& $\dim\mathcal{B}^2$ & $\dim H^2$
\\
\hline
0 & 1 & 1 & 1 & 0 \\
1 & 4 & 4 & 4 & 0 \\
2 & 6 & 5 & 5 & 0 \\
3 & 6 & 4 & 4 & 0 \\
4 & 5 & 3 & 1 & 2 \\
5 & 2 & 0 & 0 & 0 \\
\hline
\end{tabular}
\end{center}

\end{Example}

To conclude the presentation, 
in the next table, we present the speediness of the algorithm
for our two Examples~\ref{MS-Example} and~\ref{AMS-Example}, 
and also for $H^k(\frak {gl}(3),\frak{sl}(3))$:

\begin{figure}[H]
\centering {\hspace*{-0.75cm}\small
\begin{tabular}{|c||c|c|c|c|c|c|c|}
\cline{1-8} {\rm Cohomology} & Order& time(sec.) & memory(M) &
$\dim (\mathcal C^k)$ & $\dim(\mathcal Z^k)$ & $\dim(\mathcal B^k)$
& $\dim (H^k)$\\
\cline{1-8} 
Example~\ref{MS-Example}
& 2 & 0.125 & 3.6 & 70 &25&33&8 \\
\cline{1-8} 
Example~\ref{MS-Example}
& 3 & 0.125 & 4.3 & 70 &37&45&8 \\
\cline{1-8} 
Example~\ref{MS-Example}
& 4 & 0.03 & 1.4 & 35 &25&30&5 \\
\cline{1-8} 
Example~\ref{MS-Example}
& 5 & 0.0 & 0.16 & 7 &5&7&2 \\
\cline{1-8} 
Example~\ref{AMS-Example}
& 2 & 0.015 & 0.7 & 24 &15&17&2 \\
\cline{1-8} 
Example~\ref{AMS-Example}
& 3 & 0.0 & 0.18 & 8 &7&8&1 \\
\cline{1-8} $(\frak {gl}(3),\frak{sl}(3))$& 2 
& 2 & 8.6 & 252 &64&64&0 \\
\cline{1-8} $(\frak {gl}(3),\frak{sl}(3))$& 3 
& 24 & 40 & 504 &188&189&1 \\\cline{1-8}
\end{tabular}}
\end{figure}

\subsection*{Acknowledgments} 
We have the pleasure to thank Dr. Amir Hashemi for helpful discussions
during the preparation of this article.

\end{document}